\documentclass[reqno,12pt,a4paper]{amsart}

\usepackage{mathrsfs}
\usepackage{amssymb}
\usepackage{amsfonts}
\usepackage{latexsym}
\usepackage{amsthm}
\usepackage{graphicx}

\def\lb{\label}

\newcommand{\er}[1]{\textrm{(\ref{#1})}}

\begin{document}

%%%%%%%%%% Some definitions %%%%%%%%%%

%%%%%%%% Equations, theorems %%%%%%%%%
\renewcommand{\theequation}{\arabic{section}.\arabic{equation}}
\theoremstyle{plain}
\newtheorem{theorem}{\bf Theorem}[section]
\newtheorem{lemma}[theorem]{\bf Lemma}
\newtheorem{corollary}[theorem]{\bf Corollary}
\newtheorem{proposition}[theorem]{\bf Proposition}
\newtheorem{definition}[theorem]{\bf Definition}

\newtheorem{remark}[theorem]{\bf Remark}

%%%%% Alphabet %%%%%
\def\a{\alpha}  \def\cA{{\mathcal A}}     \def\bA{{\bf A}}  \def\mA{{\mathscr A}}
\def\b{\beta}   \def\cB{{\mathcal B}}     \def\bB{{\bf B}}  \def\mB{{\mathscr B}}
\def\g{\gamma}  \def\cC{{\mathcal C}}     \def\bC{{\bf C}}  \def\mC{{\mathscr C}}
\def\G{\Gamma}  \def\cD{{\mathcal D}}     \def\bD{{\bf D}}  \def\mD{{\mathscr D}}
\def\d{\delta}  \def\cE{{\mathcal E}}     \def\bE{{\bf E}}  \def\mE{{\mathscr E}}
\def\D{\Delta}  \def\cF{{\mathcal F}}     \def\bF{{\bf F}}  \def\mF{{\mathscr F}}
\def\c{\chi}    \def\cG{{\mathcal G}}     \def\bG{{\bf G}}  \def\mG{{\mathscr G}}
\def\z{\zeta}   \def\cH{{\mathcal H}}     \def\bH{{\bf H}}  \def\mH{{\mathscr H}}
\def\e{\eta}    \def\cI{{\mathcal I}}     \def\bI{{\bf I}}  \def\mI{{\mathscr I}}
\def\p{\psi}    \def\cJ{{\mathcal J}}     \def\bJ{{\bf J}}  \def\mJ{{\mathscr J}}
\def\vT{\Theta} \def\cK{{\mathcal K}}     \def\bK{{\bf K}}  \def\mK{{\mathscr K}}
\def\k{\kappa}  \def\cL{{\mathcal L}}     \def\bL{{\bf L}}  \def\mL{{\mathscr L}}
\def\l{\lambda} \def\cM{{\mathcal M}}     \def\bM{{\bf M}}  \def\mM{{\mathscr M}}
\def\L{\Lambda} \def\cN{{\mathcal N}}     \def\bN{{\bf N}}  \def\mN{{\mathscr N}}
\def\m{\mu}     \def\cO{{\mathcal O}}     \def\bO{{\bf O}}  \def\mO{{\mathscr O}}
\def\n{\nu}     \def\cP{{\mathcal P}}     \def\bP{{\bf P}}  \def\mP{{\mathscr P}}
\def\r{\rho}    \def\cQ{{\mathcal Q}}     \def\bQ{{\bf Q}}  \def\mQ{{\mathscr Q}}
\def\s{\sigma}  \def\cR{{\mathcal R}}     \def\bR{{\bf R}}  \def\mR{{\mathscr R}}
\def\S{\Sigma}  \def\cS{{\mathcal S}}     \def\bS{{\bf S}}  \def\mS{{\mathscr S}}
\def\t{\tau}    \def\cT{{\mathcal T}}     \def\bT{{\bf T}}  \def\mT{{\mathscr T}}
\def\f{\phi}    \def\cU{{\mathcal U}}     \def\bU{{\bf U}}  \def\mU{{\mathscr U}}
\def\F{\Phi}    \def\cV{{\mathcal V}}     \def\bV{{\bf V}}  \def\mV{{\mathscr V}}
\def\P{\Psi}    \def\cW{{\mathcal W}}     \def\bW{{\bf W}}  \def\mW{{\mathscr W}}
\def\o{\omega}  \def\cX{{\mathcal X}}     \def\bX{{\bf X}}  \def\mX{{\mathscr X}}
\def\x{\xi}     \def\cY{{\mathcal Y}}     \def\bY{{\bf Y}}  \def\mY{{\mathscr Y}}
\def\X{\Xi}     \def\cZ{{\mathcal Z}}     \def\bZ{{\bf Z}}  \def\mZ{{\mathscr Z}}
%*********************
\def\be{{\bf e}}
\def\bv{{\bf v}} \def\bu{{\bf u}}
\def\Om{\Omega}
%************************
\def\bbD{\pmb \Delta}
\def\mm{\mathrm m}
\def\mn{\mathrm n}
%*************************

\newcommand{\mc}{\mathscr {c}}

\newcommand{\gA}{\mathfrak{A}}          \newcommand{\ga}{\mathfrak{a}}
\newcommand{\gB}{\mathfrak{B}}          \newcommand{\gb}{\mathfrak{b}}
\newcommand{\gC}{\mathfrak{C}}          \newcommand{\gc}{\mathfrak{c}}
\newcommand{\gD}{\mathfrak{D}}          \newcommand{\gd}{\mathfrak{d}}
\newcommand{\gE}{\mathfrak{E}}
\newcommand{\gF}{\mathfrak{F}}           \newcommand{\gf}{\mathfrak{f}}
\newcommand{\gG}{\mathfrak{G}}           %\newcommand{\gg}{\mathfrak{g}}
\newcommand{\gH}{\mathfrak{H}}           \newcommand{\gh}{\mathfrak{h}}
\newcommand{\gI}{\mathfrak{I}}           \newcommand{\gi}{\mathfrak{i}}
\newcommand{\gJ}{\mathfrak{J}}           \newcommand{\gj}{\mathfrak{j}}
\newcommand{\gK}{\mathfrak{K}}            \newcommand{\gk}{\mathfrak{k}}
\newcommand{\gL}{\mathfrak{L}}            \newcommand{\gl}{\mathfrak{l}}
\newcommand{\gM}{\mathfrak{M}}            \newcommand{\gm}{\mathfrak{m}}
\newcommand{\gN}{\mathfrak{N}}            \newcommand{\gn}{\mathfrak{n}}
\newcommand{\gO}{\mathfrak{O}}
\newcommand{\gP}{\mathfrak{P}}             \newcommand{\gp}{\mathfrak{p}}
\newcommand{\gQ}{\mathfrak{Q}}             \newcommand{\gq}{\mathfrak{q}}
\newcommand{\gR}{\mathfrak{R}}             \newcommand{\gr}{\mathfrak{r}}
\newcommand{\gS}{\mathfrak{S}}              \newcommand{\gs}{\mathfrak{s}}
\newcommand{\gT}{\mathfrak{T}}             \newcommand{\gt}{\mathfrak{t}}
\newcommand{\gU}{\mathfrak{U}}             \newcommand{\gu}{\mathfrak{u}}
\newcommand{\gV}{\mathfrak{V}}             \newcommand{\gv}{\mathfrak{v}}
\newcommand{\gW}{\mathfrak{W}}             \newcommand{\gw}{\mathfrak{w}}
\newcommand{\gX}{\mathfrak{X}}               \newcommand{\gx}{\mathfrak{x}}
\newcommand{\gY}{\mathfrak{Y}}              \newcommand{\gy}{\mathfrak{y}}
\newcommand{\gZ}{\mathfrak{Z}}             \newcommand{\gz}{\mathfrak{z}}

\def\ve{\varepsilon}   \def\vt{\vartheta}    \def\vp{\varphi}    \def\vk{\varkappa}

\def\A{{\mathbb A}} \def\B{{\mathbb B}} \def\C{{\mathbb C}}
\def\dD{{\mathbb D}} \def\E{{\mathbb E}} \def\dF{{\mathbb F}} \def\dG{{\mathbb G}} \def\H{{\mathbb H}}\def\I{{\mathbb I}} \def\J{{\mathbb J}} \def\K{{\mathbb K}} \def\dL{{\mathbb L}}\def\M{{\mathbb M}} \def\N{{\mathbb N}} \def\O{{\mathbb O}} \def\dP{{\mathbb P}} \def\R{{\mathbb R}}\def\S{{\mathbb S}} \def\T{{\mathbb T}} \def\U{{\mathbb U}} \def\V{{\mathbb V}}\def\W{{\mathbb W}} \def\X{{\mathbb X}} \def\Y{{\mathbb Y}} \def\Z{{\mathbb Z}}

%%%%% Arrows %%%%%

\def\la{\leftarrow}              \def\ra{\rightarrow}            \def\Ra{\Rightarrow}
\def\ua{\uparrow}                \def\da{\downarrow}
\def\lra{\leftrightarrow}        \def\Lra{\Leftrightarrow}

%%%%% Typography %%%%%

\def\lt{\biggl}                  \def\rt{\biggr}
\def\ol{\overline}               \def\wt{\widetilde}
\def\ul{\underline}
\def\no{\noindent}

%%%%% Math signs %%%%%

\let\ge\geqslant                 \let\le\leqslant
\def\lan{\langle}                \def\ran{\rangle}
\def\/{\over}                    \def\iy{\infty}
\def\sm{\setminus}               \def\es{\emptyset}
\def\ss{\subset}                 \def\ts{\times}
\def\pa{\partial}                \def\os{\oplus}
\def\om{\ominus}                 \def\ev{\equiv}
\def\iint{\int\!\!\!\int}        \def\iintt{\mathop{\int\!\!\int\!\!\dots\!\!\int}\limits}
\def\el2{\ell^{\,2}}             \def\1{1\!\!1}
\def\sh{\sharp}
\def\wh{\widehat}
\def\bs{\backslash}
\def\intl{\int\limits}
%%%%% Math operations %%%%%

\def\na{\mathop{\mathrm{\nabla}}\nolimits}
\def\sh{\mathop{\mathrm{sh}}\nolimits}
\def\ch{\mathop{\mathrm{ch}}\nolimits}
\def\where{\mathop{\mathrm{where}}\nolimits}
\def\all{\mathop{\mathrm{all}}\nolimits}
\def\as{\mathop{\mathrm{as}}\nolimits}
\def\Area{\mathop{\mathrm{Area}}\nolimits}
\def\arg{\mathop{\mathrm{arg}}\nolimits}
\def\const{\mathop{\mathrm{const}}\nolimits}
\def\det{\mathop{\mathrm{det}}\nolimits}
\def\diag{\mathop{\mathrm{diag}}\nolimits}
\def\diam{\mathop{\mathrm{diam}}\nolimits}
\def\dim{\mathop{\mathrm{dim}}\nolimits}
\def\dist{\mathop{\mathrm{dist}}\nolimits}
\def\Im{\mathop{\mathrm{Im}}\nolimits}
\def\Iso{\mathop{\mathrm{Iso}}\nolimits}
\def\Ker{\mathop{\mathrm{Ker}}\nolimits}
\def\Lip{\mathop{\mathrm{Lip}}\nolimits}
\def\rank{\mathop{\mathrm{rank}}\limits}
\def\Ran{\mathop{\mathrm{Ran}}\nolimits}
\def\Re{\mathop{\mathrm{Re}}\nolimits}
\def\Res{\mathop{\mathrm{Res}}\nolimits}
\def\res{\mathop{\mathrm{res}}\limits}
\def\sign{\mathop{\mathrm{sign}}\nolimits}
\def\span{\mathop{\mathrm{span}}\nolimits}
\def\supp{\mathop{\mathrm{supp}}\nolimits}
\def\Tr{\mathop{\mathrm{Tr}}\nolimits}
\def\BBox{\hspace{1mm}\vrule height6pt width5.5pt depth0pt \hspace{6pt}}

%%%%%%%%%%%%% specialities %%%%%%%%%%%%%%

\newcommand\nh[2]{\widehat{#1}\vphantom{#1}^{(#2)}}
%{{\mathop{#1}\limits^\wedge}\vphantom{#1}^{(#2)}}
\def\dia{\diamond}

\def\Oplus{\bigoplus\nolimits}

%%%%%%%%%%% End of definitions %%%%%%%%%%

%%%%% OLD OLD OLD

\def\qqq{\qquad}
\def\qq{\quad}
\let\ge\geqslant
\let\le\leqslant
\let\geq\geqslant
\let\leq\leqslant
\newcommand{\ca}{\begin{cases}}
\newcommand{\ac}{\end{cases}}
\newcommand{\ma}{\begin{pmatrix}}
\newcommand{\am}{\end{pmatrix}}
\renewcommand{\[}{\begin{equation}}
\renewcommand{\]}{\end{equation}}
\def\eq{\begin{equation}}
\def\qe{\end{equation}}
\def\[{\begin{equation}}
\def\bu{\bullet}

\title[Invariants of magnetic Laplacians on periodic graphs]
{Invariants of magnetic Laplacians on periodic graphs}

\date{\today}
\author[Evgeny Korotyaev]{Evgeny Korotyaev}
\address{Dep. of Mathematical Analysis, Saint-Petersburg State University, Universitetskaya nab. 7/9, St. Petersburg, 199034, Russia,
\ korotyaev@gmail.com, \
e.korotyaev@spbu.ru,}
\author[Natalia Saburova]{Natalia Saburova}
\address{Dep. of Mathematical Analysis, Algebra and Geometry, Northern (Arctic) Federal University, Severnaya Dvina Emb. 17, Arkhangelsk, 163002, Russia,
 \ n.saburova@gmail.com, \ n.saburova@narfu.ru}

\subjclass{} \keywords{discrete magnetic Laplacians, periodic graphs, spectral bands}

\begin{abstract}
We consider a magnetic Laplacian with periodic magnetic potentials on periodic discrete graphs. Its spectrum consists of a finite number of bands, where degenerate bands are eigenvalues of infinite multiplicity. We obtain a specific decomposition of the magnetic Laplacian into a direct integral in terms of minimal forms. A minimal form is a periodic function defined on edges of the periodic graph with a minimal support on the period. It is crucial that fiber magnetic Laplacians (matrices) have the minimal number of coefficients depending on the quasimomentum and the minimal number of coefficients depending on the magnetic potential. We show that these numbers are invariants for the magnetic Laplacians on periodic graphs. Using this decomposition, we estimate the position of each band, the Lebesgue measure of the magnetic Laplacian spectrum and a variation of the spectrum under a perturbation by a magnetic field in terms of these invariants and minimal forms. In addition, we consider an inverse problem: we determine necessary and sufficient conditions for matrices depending on the quasimomentum on a finite graph to be fiber magnetic Laplacians. Moreover, similar results for magnetic Schr\"odinger operators with periodic potentials are obtained.
\end{abstract}

\maketitle

\section {\lb{Sec1}Introduction}
\setcounter{equation}{0}
Laplace and Schr\"odinger operators on periodic discrete graphs are of interest due to their applications to problems of physics and chemistry.
We consider Laplacians and Schr\"odinger operators on periodic discrete graphs in external magnetic fields. The magnetic Laplacian is expressed in terms of a magnetic vector potential which is a function of the graph edges. We assume that magnetic and electric potentials are periodic. This guarantees a finite band structure of the spectrum. It is well known that in this case the
spectrum of the magnetic operators consists of an absolutely continuous part
(a~union of a finite number of non-degenerate bands) and a finite
number of flat bands, i.e., eigenvalues  of infinite multiplicity.

A discrete analogue of the magnetic Laplacian on $\R^2$ was originally introduced by Harper \cite{H55}. It describes the behavior of an electron moving on the square lattice $\Z^2$ in an external uniform magnetic field perpendicular to the lattice in the so-called tight-binding approximation. This magnetic Laplacian on $\Z^2$ admits a reduction to the Harper operator. The graphical presentation of the dependence of the spectrum of this operator on the magnetic flux through the unit cell of the lattice is known as the Hofstadter butterfly \cite{Ho76}.  Helffer and Sj\"ostrand in a series of papers \cite{HS88}, \cite{HS89}, \cite{HS90} described some spectral properties of the Harper operator. An algebraic approach to this operator was put forward by Bellissard (for more details see \cite{Be92}, \cite{Be94}). Note that there are results about the Hofstadter-type spectrum of the Harper model on other planar graphs (the triangular, hexagonal, Kagome lattices) (see \cite{Hou09}, \cite{Ke92}, \cite{KeR14}, \cite{HKeR16}, and the references therein).

Discrete magnetic Laplacians on graphs were studied by many authors (see, e.g., \cite{LL93}, \cite{S94} and the references therein). Magnetic Schr\"odinger operators on periodic discrete graphs were considered in \cite{HS99a}, \cite{KS17}. Higuchi and Shirai \cite{HS99a} studied the behavior of the bottom of the spectrum as a function of the magnetic flux. Korotyaev and Saburova \cite{KS17} estimated the Lebesgue measure of the spectrum in terms of Betti numbers defined by \er{benu}. They also obtained estimates of a variation of the spectrum of the Schr\"odinger operators under a perturbation by a magnetic field in terms of magnetic fluxes and estimates of effective masses associated with the ends of each spectral band for magnetic Laplacians in terms of geometric parameters of the graphs.

In the present paper we consider magnetic Schr\"odinger operators with periodic magnetic and electric potentials on periodic discrete graphs. We will essentially use a notion of \emph{minimal forms} on graphs introduced in \cite{KS18}. A minimal form is a periodic function defined on edges of the periodic graph with a minimal support on the period. We describe our main goals:

$\bu$  to decompose the magnetic Schr\"odinger operators into a direct integral, where fiber operators are expressed in terms of minimal forms and have the minimal number $2\cI$ of coefficients depending on the quasimomentum and the minimal number of coefficients depending on the magnetic potential;

$\bu$ to show that these numbers are invariants for the magnetic Schr\"odinger operators on periodic graphs and the difference between the Betti number and each of these invariants can be arbitrarily large (for specific graphs);

$\bu$ to determine a localization of bands in terms of eigenvalues of the magnetic Schr\"odinger operator on a finite graph obtained from the periodic one by deleting a support of a minimal form;

$\bu$ to obtain an estimate of the Lebesgue measure of the spectrum in terms of the invariant $\cI$; to show that this estimate becomes an identity for specific graphs;

$\bu$ to estimate the variation of the spectrum of the Schr\"odinger operators under a perturbation by a magnetic field in terms of minimal forms;

$\bu$ to solve an inverse problem: to determine necessary and sufficient conditions for matrices depending on the quasimomentum on a finite graph to be fiber magnetic Laplacians.

\subsection{Magnetic Schr\"odinger operators on periodic graphs.}
Let $\cG=(\cV,\cE)$ be a connected infinite graph, possibly  having
loops and multiple edges and embedded into the space $\R^d$. Here $\cV$ is the set of its vertices and $\cE$ is the set of its unoriented edges. Considering each edge in $\cE$ to have two orientations, we introduce the set $\cA$ of all oriented edges.
An edge starting at a vertex $u$ and ending at a vertex $v$ from $\cV$ will be denoted as the ordered pair $(u,v)\in\cA$ and is said to be \emph{incident} to the vertices. Let $\ul\be=(v,u)$ be the inverse edge of $\be=(u,v)\in\cA$. Vertices $u,v\in\cV$ will be called \emph{adjacent} and denoted by $u\sim v$, if $(u,v)\in \cA$. We define the degree ${\vk}_v$ of
the vertex $v\in\cV$ as the number of all edges from
$\cA$, starting at $v$. A sequence of directed edges $(\be_1,\be_2,\ldots,\be_n)$ is called a \emph{path} if the terminus of the edge $\be_s$ coincides with the origin of the edge $\be_{s+1}$ for all $s=1,\ldots,n-1$. If the terminus of $\be_n$ coincides with the origin of $\be_1$, the path is called a \emph{cycle}.

Let $\G$ be a lattice of rank $d$ in $\R^d$ with a basis $\A=\{a_1,\ldots,a_d\}$, i.e.,
$$
\G=\Big\{a : a=\sum_{s=1}^dn_sa_s, \; n_s\in\Z,\; s\in\N_d\Big\}, \qqq \N_d=\{1,\ldots,d\},
$$
and let
\[\lb{fuce}
\Omega=\Big\{x\in\R^d : x=\sum_{s=1}^dt_sa_s, \; 0\leq t_s<1,\; s\in\N_d\Big\}
\]
be the \emph{fundamental cell} of the lattice $\G$. We define the equivalence relation on $\R^d$:
$$
x\equiv y \; (\hspace{-4mm}\mod \G) \qq\Leftrightarrow\qq x-y\in\G \qqq
\forall\, x,y\in\R^d.
$$

We consider \emph{locally finite $\G$-periodic graphs} $\cG$, i.e., graphs satisfying the
following conditions:
\begin{itemize}
  \item[1)] $\cG=\cG+a$ for any $a\in\G$;
  \item[2)] the quotient graph  $G_*=\cG/\G$ is finite.
\end{itemize}
The basis $a_1,\ldots,a_d$ of the lattice $\G$ is called the {\it periods} of $\cG$. We also call the quotient graph $G_*=\cG/\G$ the \emph{fundamental graph}
of the periodic graph $\cG$. The fundamental graph $G_*$ is a graph on the $d$-dimensional torus $\R^d/\G$. The graph $G_*=(\cV_*,\cE_*)$ has the vertex set $\cV_*=\cV/\G$, the set $\cE_*=\cE/\G$ of unoriented edges and the set $\cA_*=\cA/\G$ of oriented edges which are finite.

Let $\ell^2(\cV)$ be the Hilbert space of all square summable functions $f:\cV\to \C$ equipped with the norm
$$
\|f\|^2_{\ell^2(\cV)}=\sum_{v\in\cV}|f(v)|^2<\infty.
$$

We define the so-called \emph{combinatorial magnetic Laplacian} $\D_\a$ on
$f\in\ell^2(\cV)$ by
\[
\lb{DLO}
\big(\D_{\a} f\big)(v)=\sum_{\be=(v,u)\in\cA}\big(f(v)-e^{i\a(\be)}f(u)\big), \qqq v\in\cV,
\]
where $\a:\cA\ra\R$ is a $\G$-periodic \emph{magnetic vector potential on $\cG$}, i.e., it satisfies for all $(\be,a)\in \cA\ts\G$:
\[\lb{MPot}
\a(\ul\be\,)=-\a(\be),\qqq \a(\be+a)=\a(\be).
\]
The sum in \er{DLO} is taken over all oriented edges starting at the vertex $v$.
If $\a=0$, then $\D_0$ is just the usual Laplacian $\D$ without magnetic potentials:
\[\lb{DLOW}
\big(\D f\big)(v)=\sum_{(v,u)\in\cA}\big(f(v)-f(u)\big), \qqq v\in\cV.
\]

\medskip

It is well known (see, e.g., \cite{HS99a}) that $\D_{\a}$
is self-adjoint and its spectrum $\s(\D_{\a})$ is contained in $[0,2\vk_+]$, i.e.,
\[
\lb{bf}
\begin{aligned}
\s(\D_{\a})\subset[0,2\vk_+],\qqq
\textrm{where}\qqq
\vk_+=\sup_{v\in\cV}\vk_v<\iy.
\end{aligned}
\]

\no \textbf{Remark.} There are other definitions of discrete magnetic Laplacians on graphs: weighted, normalized, standard Laplacians. For more details see, e.g., \cite{KS17} and references therein.

\medskip

We consider the magnetic Schr\"odinger operator $H_{\a}$ acting on the Hilbert space $\ell^2(\cV)$ and given by
\[
\lb{Sh}
H_{\a}=\D_{\a}+Q,
\]
where $Q$ is a real $\G$-periodic potential, i.e., it satisfies
for all $(v,a)\in\cV\ts\G$:
\[
\lb{Pot} \big(Q f\big)(v)=Q(v)f(v),\qqq Q(v+a)=Q(v).
\]

\subsection{Minimal forms} A notion of minimal forms on graphs was introduced in \cite{KS18}. The authors used minimal forms on the fundamental graphs in order to obtain some spectral estimates for Schr\"odinger operators without a magnetic field. For convenience of the reader we repeat here the definition of the minimal forms on graphs, since they will be used in the formulation of our results.

Let $G=(V,E)$ be a finite connected graph with the vertex set $V$, the set $E$ of unoriented edges and the set $A$ of oriented edges. Denote by $\cC$ the cycle space of the graph $G$.

A vector-valued function $\gb: A\ra\R^n$, $n\in\N$, satisfying the condition $\gb(\ul\be\,)=-\gb(\be)$ for all $\be\in A$ is called a \emph{1-form} on the graph $G$. For any 1-form $\gb$ we define the vector-valued \emph{flux function} $\Phi_\gb:\cC\to\R^n$ as follows:
\[\lb{mafla}
\Phi_\gb(\mathbf{c})=\sum\limits_{\be\in\mathbf{c}}\gb(\be),  \qqq \mathbf{c}\in\cC.
\]
We call $\Phi_\gb(\mathbf{c})\in\R^n$ the \emph{flux} of the 1-form $\gb$ through the cycle $\mathbf{c}$.

We fix a 1-form $\mathbf{x}:A\ra\R^n$ and define the set $\mF(\mathbf{x})$ of all 1-forms $\gb:A\ra\R^n$ satisfying the condition $\Phi_\gb=\Phi_{\mathbf{x}}$:
\[\lb{vv1f}
\mF(\mathbf{x})=\{\gb: A\ra\R^n : \gb(\ul\be\,)=-\gb(\be)\qq \textrm{for all} \qq \be\in A, \qq \textrm{and} \qq \Phi_\gb=\Phi_{\mathbf{x}}\},
\]
where $\Phi_{\mathbf{x}}$ is the flux function of the form $\mathbf{x}$ defined by \er{mafla}.

In the set $\mF(\mathbf{x})$ we specify a subset of \emph{minimal} 1-forms. Let $\#\cM$ denote the number of elements in a set $\cM$. A 1-form $\gm\in\mF(\mathbf{x})$ is called \emph{minimal} if
\[\lb{dmf}
\#\supp \gm\leq\#\supp \gb \qqq \textrm{for any 1-form} \qqq \gb\in\mF(\mathbf{x}).
\]
Among all 1-forms $\gb\in\mF(\mathbf{x})$ the minimal form $\gm$ has a minimal support. Since any 1-form $\gb\in\mF(\mathbf{x})$ has a finite support, minimal forms $\gm\in\mF(\mathbf{x})$ exist. An explicit expression for all minimal forms $\gm\in\mF(\mathbf{x})$ and for the number of edges in their supports are given in Theorem \ref{Pphi}.

\subsection{Magnetic and coordinate forms} For each $x\in\R^d$ we introduce the vector $x_\A\in\R^d$ by
\[\lb{cola}
x_\A=(t_1,\ldots,t_d), \qqq \textrm{where} \qq x=\textstyle\sum\limits_{s=1}^dt_sa_s.
\]
In other words, $x_\A$ is the coordinate vector of $x$ with respect to the basis $\A=\{a_1,\ldots,a_d\}$ of the lattice $\G$.

For any oriented edge $\be=(u,v)\in\cA$ of the periodic graph $\cG$ the \emph{edge coordinates} $\k(\be)$ are defined as the vector in $\R^d$ given by
\[
\lb{edco}
\k(\be)=v_\A-u_\A\in\R^d.
\]

The edge coordinates $\k(\be)$ and the magnetic potential $\a$ are $\G$-periodic, i.e., they satisfy
\[\lb{Gpe}
\k(\be+a)=\k(\be),\qqq \a(\be+a)=\a(\be),\qqq \forall\, (\be,a)\in\cA \ts\G.
\]
On the set $\cA$ of all oriented edges of the $\G$-periodic graph $\cG$ we define the surjection
\[\lb{sur}
\gf:\cA\rightarrow\cA_*=\cA/\G,
\]
which maps each $\be\in\cA$ to its equivalence class $\be_*=\gf(\be)$ which is an oriented edge of the fundamental graph $G_*$.

The identities \er{Gpe} allow us to define uniquely two 1-forms on the fundamental graph $G_*$:

$\bu$ the vector-valued \emph{coordinate form} $\k: \cA_*\ra\R^d$ is induced by the
coordinates of the periodic graph edges:
\[
\lb{dco}
\k(\be_*)=\k(\be) \qq \textrm{ for some $\be\in\cA$ \; such that }  \; \be_*=\gf(\be), \qqq \be_*\in\cA_*;
\]

$\bu$ the scalar-valued \emph{magnetic form} $\a:\cA_*\ra\R$ is induced by the
magnetic potential $\a$:
\[
\lb{Mdco}
\a(\be_*)=\a(\be)\qq \textrm{ for some $\be\in\cA$ \; such that }  \; \be_*=\gf(\be), \qqq \be_*\in\cA_*,
\]
where $\gf$ is defined by \er{sur}.

Since the magnetic potential $\a$ appears in the Laplacian $\D_\a$ via the factor $e^{i\a(\be)}$, $\be\in\cA$, we will consider fluxes of the magnetic form $\a$ modulo $2\pi$, i.e., for the magnetic form $\a:\cA_*\to\R$ we define the flux function $\Phi_\a:\cC\to(-\pi,\pi]$ by
\[\lb{modpi}
\Phi_\a(\mathbf{c})=\bigg(\sum\limits_{\be\in\mathbf{c}}\a(\be)\bigg)\hspace{-3mm} \mod 2\pi, \qqq \Phi_\a(\mathbf{c})\in(-\pi,\pi], \qqq \mathbf{c}\in\cC,
\]
where $\cC$ is the cycle space of the fundamental graph $G_*$.

\section {\lb{Sec1'} Main results}
\setcounter{equation}{0}
\subsection{Decomposition of magnetic Schr\"odinger operators}
Recall that in the standard case the magnetic Laplacian $\D_A=\big(-i\na-A(x)\big)^2$ in $L^2(\R^d)$, where $A(x)$ is a $\G$-periodic magnetic potential, $x\in\R^d$, is unitarily equivalent to the constant fiber direct integral:
$$
U\D_AU^{-1}=\int^\oplus_{\T^d}\D_A(\vt)d\vt, \qqq
\D_A(\vt)=\big(-i\na-A(x)+\vt\big)^2,\qqq \T^d=\R^d/(2\pi\Z)^d.
$$
The parameter $\vt$ is called the \emph{quasimomentum}, $U:L^2(\R^d)\to \int^\oplus_{\T^d}L^2(\Omega)d\vt$ is the Gelfand transformation, and $\Omega$ is the fundamental cell of the lattice $\G$. The fiber operator $\D_A(\vt)$ is a magnetic Laplacian with the magnetic potential $-A(x)+\vt$, on $L^2(\Omega)$.

In the case of Laplace and Schr\"odinger operators on periodic graphs we have the similar decomposition (see Theorem \ref{TFD1}). Here we present a specific direct integral for magnetic Schr\"odinger operators on periodic graphs, where fiber operators have the minimal number of coefficients depending on the quasimomentum and the minimal number of coefficients depending on the magnetic potential.

We introduce the Hilbert space
\[\lb{Hisp}
\mH=L^2\Big(\T^{d},{d\vt\/(2\pi)^d}\,,\cH\Big)=\int_{\T^{d}}^{\os}\cH\,{d\vt
\/(2\pi)^d}\,, \qqq \cH=\ell^2(\cV_*),
\]
i.e., a constant fiber direct integral equipped with the norm
$$
\|g\|^2_{\mH}=\int_{\T^{d}}\|g(\vt,\cdot)\|_{\ell^2(\cV_*)}^2\frac{d\vt
}{(2\pi)^d}\,,
$$
where the function $g(\vt,\cdot)\in\cH$ for almost all
$\vt\in\T^{d}$.

\begin{theorem}
\label{TDImf}
Let $(\gm,\phi)\in\mF(\k)\ts\mF(\a)$ be two minimal forms on the fundamental graph $G_*=(\cV_*,\cE_*)$, where $\mF(\mathbf{x})$ is given by \er{vv1f}; $\k$ and $\a$ are the coordinate and magnetic forms defined by \er{edco}, \er{dco} and \er{Mdco}, respectively. Then

i) The image of the minimal form $\gm$ generates the group $\Z^d$. The image of its flux function $\F_\gm$ defined by \er{mafla} satisfies $\F_\gm(\cC)=\Z^d$, where
$\cC$ is the cycle space of $G_*$.

ii) The magnetic Schr\"odinger operator $H_{\a}=\D_{\a}+Q$ on $\ell^2(\cV)$ has the following decomposition into a constant fiber direct integral
\[
\lb{raz1m}
\mU_{\gm,\f} H_{\a} \mU_{\gm,\f}^{-1}=\int^\oplus_{\T^{d}}H_{\gm,\f}(\vt)\frac{d\vt
}{(2\pi)^d}\,,
\]
where the unitary operator $\mU_{\gm,\f}:\ell^2(\cV)\to\mH$ is a product of the Gelfand transformation $U$ and the gauge transformation $W_{\gm,\f}$ defined by \er{5001} and \er{Uvt}, respectively.  Here the fiber magnetic Schr\"odinger operator $H_{\gm,\f}(\vt)$ and the fiber magnetic Laplacian $\D_{\gm,\f}(\vt)$ are given by
\[
\label{Hvtm}
H_{\gm,\f}(\vt)=\D_{\gm,\f}(\vt)+Q,\qqq \forall\,\vt\in\T^{d},
\]
\[
\label{l2.13am}
\begin{aligned}
\big(\D_{\gm,\f}(\vt)f\big)(v)=\vk_vf(v)
-\sum_{\be=(v,u)\in\cA_*}e^{i(\phi(\be)+\lan\gm(\be),\vt\ran)}f(u), \qq f\in\ell^2(\cV_*), \qq v\in\cV_*,
\end{aligned}
\]
where $\lan\cdot\,,\cdot\ran$ denotes the standard inner product in $\R^d$.

iii) Let $(\gm',\phi')\in\mF(\k)\ts\mF(\a)$ be two other minimal forms on $G_*$. Then for each $\vt\in\T^{d}$ the operators $\D_{\gm,\f}(\vt)$ and $\D_{\gm',\f'}(\vt)$ are unitarily equivalent.

iv) In the identity \er{l2.13am} for the fiber Laplacian $\D_{\gm,\f}(\cdot)$

$\bu$ the number of exponents $e^{i\lan\gm(\be),\,\cdot\,\ran}\neq1$, $\be\in\cA_*$, is equal to $2\cI$;

$\bu$ the number of exponents $e^{i\phi(\be)}\neq1$, $\be\in\cA_*$, is equal to $2\cI_\a$;

$\bu$ the number of exponents $e^{i(\phi(\be)+\lan\gm(\be),\,\cdot\,\ran)}\neq1$, $\be\in\cA_*$, is equal to $2\cI_{\gm,\f}$,
where
\[\lb{dIm}
\cI=\textstyle\frac12\,\#\supp \gm,\qqq  \cI_\a=\textstyle\frac12\,\#\supp\f,
\]
\[\lb{decI}
\cI_{\gm,\f}=\textstyle\frac12\,\#(\supp \gm\cup\supp\phi).
\]

v) The number $\cI$ \textbf{does not depend} on any of the following choices

$\bu$ the choice of the minimal form $\gm\in\mF(\k)$;

$\bu$ the choice of the embedding of $\cG$ into $\R^d$;

$\bu$ the choice of the basis $a_1,\ldots,a_d$ of the lattice $\G$.
\end{theorem}

\no \textbf{Remarks.} 1) For simple periodic graphs (the $d$-dimensional lattice, the hexagonal lattice, the Kagome lattice, etc.) it is not difficult to find  minimal forms $\gm\in\mF(\k)$ and $\phi\in\mF(\a)$ on their fundamental graphs. But for an arbitrary periodic graph this may be a rather complicated problem (see Theorem \ref{Pphi}).

2) We can consider the fiber operator $\D_{\gm,\f}(\vt)$, $\vt\in
\T^{d}$, as the magnetic Laplacian given by \er{DLO} with the
magnetic vector potential $\a(\be)=\f(\be)+\lan\gm(\be),\vt\ran$ defined on edges $\be$ of the fundamental graph $G_*$.

3) A decomposition of the magnetic Schr\"odinger operator
on periodic discrete graphs into a constant fiber direct integral was given in \cite{KS17}, where the fiber operator was expressed in terms of the magnetic form $\a$ and a special 1-form $\t\in\mF(\k)$ defined by \er{in}, \er{inf} (see Theorem \ref{TFD1}). In this representation the fiber operator depends on $\#\supp\a$ values $\a(\be)$ of the magnetic form $\a$ on edges $\be\in\supp\a$ and $\#\supp\t$ values $\t(\be)$ of the form $\t$ on $\be\in\supp\t$. These numbers $\#\supp\a$ and $\#\supp\t$ are difficult to control. Moreover, the number $\#\supp\t$  essentially depends on the choice of the embedding of the periodic graph $\cG$ into the space $\R^d$ and both the numbers $\#\supp\a$ and $\#\supp\t$ can be significantly greater than, respectively, $2\cI_\a$ and $2\cI$ (see \er{dIm}), which are invariants for the magnetic Laplacian $\D_\a$ on $\cG$.

4) In Theorem \ref{TDIBB} we give a general decomposition of the magnetic Schr\"odinger operator on periodic discrete graphs into a constant fiber direct integral, where the fiber operator $\D_{\gb,\ga}(\vt)$ is expressed in terms of any 1-forms $(\gb,\ga)\in\mF(\k)\ts\mF(\a)$. Moreover, we show that among all fiber magnetic Laplacians $\D_{\gb,\ga}(\vt)$ the fiber operator $\D_{\gm,\f}(\vt)$ (which is not uniquely defined and depends on the choice of the minimal forms $(\gm,\f)\in\mF(\k)\ts\mF(\a)$) has the minimal number $2\cI$ of coefficients depending on the quasimomentum $\vt$ and the minimal number $2\cI_\a$ of coefficients depending on the form $\f$. The numbers $\cI$ and $\cI_\a$ are uniquely defined and do not depend on the choice of $\gm$ and $\f$.

5) The set $\mF(\k)$ is defined in terms of the coordinate form $\k$. But we can introduce this set using another its representative, for example, the 1-form $\t\in\mF(\k)$ defined by \er{in}, \er{inf}.

6) The minimal forms $\gm$, of course, depend on the choice of the embedding of the periodic graph $\cG$ into the space $\R^d$ and on the choice of the basis $a_1,\ldots,a_d$ of the lattice $\G$. But, due to Theorem \ref{TDImf}.\emph{v}) the number $2\cI$ of entries in $\#\supp \gm$ does not depend on these choices. Thus, $\cI$ is an \emph{invariant} of the periodic graph $\cG$. We note that the number $\cI_\a$ defined in \er{dIm} is also an \emph{invariant} of the magnetic Laplacian $\D_\a$, since it does not depend on the choice of the minimal form $\f\in\mF(\a)$.

\medskip

We formulate some simple properties of the invariants $\cI$ and $\cI_\a$ defined by \er{dIm}.

\begin{proposition} \lb{TNNI0}
i) The invariant $\cI$ of the periodic graph $\cG$ satisfies
\[\lb{nni0}
d\leq\cI\leq\b, \qqq \b=\#\cE_*-\#\cV_*+1,
\]
where $\b$ is the Betti number of the fundamental graph $G_*=(\cV_*,\cE_*)$.

ii) Let $\cG$ be a periodic graph. Then for any integer $\mn\in[0,\b]$, there exists a magnetic potential $\a$ such that the invariant $\cI_\a$ of the magnetic Laplacian $\D_\a$ on $\cG$ is equal to $\mn$.

iii) For any nonnegative integer $M$ and positive integer $N$, there exist a periodic graph $\cG$ and a magnetic potential $\a$ such that
\[\lb{nni00}
\b-\cI=M, \qqq \b-\cI_\a=N.
\]

iv) For any nonnegative integer $K$, there exist a periodic graph $\cG$ and a magnetic potential $\a$ such that
\[\lb{ni00}
\b-\cI_{\gm,\f}=K,
\]
where $\cI_{\gm,\f}$ is defined by \er{decI}.
\end{proposition}

\no \textbf{Remark.} The inequalities in \er{nni0} may become identities or strict inequalities (for some specific periodic graphs). For more details see Proposition 2.2 in \cite{KS18}.

\medskip

Now we consider the inverse problem: when a matrix-valued function $\mA(\vt)$ depending on the quasimomentum $\vt$ on a finite connected graph $G$ is a fiber magnetic Laplacian for some periodic graph $\cG$.

\begin{corollary}\label{TCo1}
Let $\gm: A\ra\R^d$, $\f: A\ra(-\pi,\pi]$ be minimal forms on a finite connected graph $G=(V,E)$ satisfying $\Phi_{\gm}(\cC)=\Z^d$, where $\cC$ is the cycle space of $G$, and $\Phi_{\gm}$ is the flux function of the form $\gm$ given by \er{mafla}. Define the operator $\mA(\vt):\ell^2(V)\to\ell^2(V)$, $\vt\in\T^d$, by
\[
\label{l2.13at}
\begin{aligned}
\big(\mA(\vt)f\big)(v)=\vk_vf(v)
-\sum_{\be=(v,u)\in A}e^{i(\phi(\be)+\lan\gm(\be),\vt\ran)}f(u), \qqq v\in V.
\end{aligned}
\]
Then $\mA(\vt)$ is a fiber operator for the magnetic Laplacian $\D_\f$ with the magnetic potential $\f$ on some periodic graph $\cG$ with the fundamental graph $G_*$ isomorphic to $G$.
\end{corollary}

\no \textbf{Remark.} This corollary and Theorem \ref{TDImf}.\emph{i,ii}) give  necessary and sufficient conditions for the operator of the form \er{l2.13at} to be a fiber magnetic Laplacian on a given graph $G$.

\medskip

In \er{l2.13am} the fiber magnetic Laplacian $\D_{\gm,\f}(\vt)$ depends on $2\cI_\a$ parameters $\f(\be)$, $\be\in\supp\f$. Now we show that using a simple change of variables we can reduce their number. Due to Theorem \ref{TDImf}.\emph{i}), there exist edges $\be_1,\ldots,\be_d\in\supp\gm$ with linear independent values $\gm(\be_1),\ldots,\gm(\be_d)$ of the form $\gm$. We denote by $\cB$ the set
\[
\lb{wulB}
\cB=\{\be_1,\ldots,\be_d,\ul\be_1,\ldots,\ul\be_d\}.
\]

\begin{corollary}\label{TCo0}  Let $(\gm,\phi)\in\mF(\k)\ts\mF(\a)$ be two minimal forms on the fundamental graph $G_*=(\cV_*,\cE_*)$, where $\mF(\mathbf{x})$ is given by \er{vv1f}, $\k$ and $\a$ are the coordinate and magnetic forms defined by \er{edco}, \er{dco} and \er{Mdco}, respectively. Then

i) The magnetic Schr\"odinger operator $H_{\a}=\D_{\a}+Q$ on $\ell^2(\cV)$ has the following decomposition into a constant fiber direct integral
\[
\lb{raz2m}
\mU_{\gm,\f} H_{\a} \mU_{\gm,\f}^{-1}=\int^\oplus_{\T^{d}}H_{\gm,\wt\f}(\vt)\frac{d\vt
}{(2\pi)^d}\,,
\]
where the unitary operator $\mU_{\gm,\f}:\ell^2(\cV)\to\mH$ is defined in Theorem \ref{TDImf}.  Here the fiber magnetic Schr\"odinger operator $H_{\gm,\wt\f}(\vt)$ and the fiber magnetic Laplacian $\D_{\gm,\wt\f}(\vt)$ are given by
\[
\label{Hvtm2}
H_{\gm,\wt\f}(\vt)=\D_{\gm,\wt\f}(\vt)+Q,\qqq \forall\,\vt\in\T^{d},
\]
\[
\label{mmfr}
\big(\D_{\gm,\wt\f}(\vt)f\big)(v)=\vk_vf(v)
-\sum_{\be=(v,u)\in\cA_*}e^{i(\wt\phi(\be)+\lan\gm(\be),\vt\ran)}f(u), \qq f\in\ell^2(\cV_*), \qq v\in \cV_*,
\]
where the 1-form $\wt\f$ is defined by
\[\lb{wtcat}
\wt\f(\be)=\left\{
\begin{array}{cl}
 \big(\phi(\be)+\lan\gm(\be),\vt_0\ran\big) , &  \textrm{ if } \, \be\in (\supp \gm\cup\supp\phi)\sm\cB\\[6pt]
  0, \qq & \textrm{ otherwise } \\
\end{array}\right.,
\]
for some $\vt_0\in\T^d$, and $\cB$ given by \er{wulB}.

ii) The number of exponents $e^{i\wt\phi(\be)}\neq1$, $\be\in\cA_*$, in the identity \er{mmfr} for the fiber Laplacian $\D_{\gm,\wt\f}(\cdot)$ is less than or equal to $2(\cI_{\gm,\f}-d)$, where $\cI_{\gm,\f}$ is defined by \er{decI}. In particular, if $\cI_{\gm,\f}=d$, then the magnetic Schr\"odinger operator $H_{\a}$ is unitarily equivalent to the Schr\"odinger operator $H_0$ without a magnetic field.
\end{corollary}

\no \textbf{Remarks.} 1) We note that this proposition holds for any 1-forms $(\gm,\phi)\in\mF(\k)\ts\mF(\a)$ (including minimal ones). In particular, in our earlier paper \cite{KS17} we proved this direct integral decomposition for a special 1-form $\t\in\mF(\k)$ defined by \er{in}, \er{inf}, and a 1-form $\a_*\in\mF(\a)$, equal to zero on edges of a spanning tree of the fundamental graph $G_*$ and coinciding with fluxes of $\a$ on other edges of $G_*$. However, $\t$ and $\a_*$, in general, are not minimal forms on $G_*$.

2) \lb{Rem2} Let $\phi\in\mF(\a)$ be a minimal form. Using Theorem \ref{Pphi}, it can be shown that there always exists a form $\gm\in\mF(\k)$ (in general, not minimal) such that $\cI_{\gm,\f}=\frac12\,\#(\supp \gm\cup\supp\phi)\leq\b$.

3) It is known (see, e.g., \cite{HS99b}) that if $\b=d$, then the magnetic Schr\"odinger operator $H_{\a}$ is unitarily equivalent to the Schr\"odinger operator $H_0$ without a magnetic field. Due to Proposition \ref{TNNI0}.\emph{iv}),  $\cI_{\gm,\f}$ can be arbitrarily smaller than $\b$. This and the previous remark yield that the sufficient condition $\cI_{\gm,\f}=d$ for the unitary equivalence of $H_\a$ and $H_0$ in Corollary \ref{TCo0}.\emph{ii}) is stronger than the condition $\b=d$. We note that this sufficient condition $\cI_{\gm,\f}=d$ can not be obtained from the decomposition in \cite{KS17}.

4) For the hexagonal lattice and the $d$-dimensional lattice with minimal fundamental graphs, $\b=d$. Then the spectrum $\s(\D_{\a})$ of the magnetic Laplacian $\D_{\a}$ on these graphs satisfies $\s(\D_{\a})=\s(\D_0)=[0,2\vk_+]$, where
$\vk_+$ is the degree of each vertex of the graph, i.e., $\s(\D_{\a})$ does not depend on the magnetic potential $\a$. If the fundamental graphs are not minimal, then, generally speaking, $\D_{\a}$ depends on $\a$.

\subsection{Spectrum of the magnetic Schr\"odinger operator.}
Theorem \ref{TDImf} and standard arguments (see Theorem XIII.85 in
\cite{RS78}) describe the spectrum of the magnetic Schr\"odinger operator
$H_{\a}=\D_{\a}+Q$. Each fiber operator $H_{\gm,\f}(\vt)$, $\vt\in\T^{d}$,
has $\n$ eigenvalues $\l_{\a,n}(\vt)$, $n\in\N_\n$, $\n=\#\cV_*$, which are labeled in non-decreasing order (counting multiplicities) by
\[
\label{eq.3} \l_{\a,1}(\vt)\leq\l_{\a,2}(\vt)\leq\ldots\leq\l_{\a,\nu}(\vt),
\qqq \forall\,\vt\in\T^{d}.
\]
Since $H_{\gm,\f}(\vt)$ is self-adjoint and analytic in $\vt\in\T^{d}$, each $\l_{\a,n}(\cdot)$, $n\in\N_\n$, is a real and piecewise analytic function on the torus $\T^{d}$ and creates the \emph{spectral band} (or \emph{band} for short) $\s_n(H_{\a})$ given by
\[
\lb{ban.1}
\s_{n}(H_{\a})=[\l_{\a,n}^-,\l_{\a,n}^+]=\l_{\a,n}(\T^{d}).
\]
\no \textbf{Remark.} From \er{ban.1} it follows that
\[\lb{HiSi}
\l_{\a,1}^{-}\leq\l_{\a,1}(0).
\]
Note that if there is no magnetic field, that is $\a=0$, Sy and Sunada \cite{SS92} proved
that $\l_{0,1}^{-}=\l_{0,1}(0)$. However, the equality in \er{HiSi} does not hold for general $\a$, since for some specific graphs we have the strict inequality (see Examples 5.2 -- 5.6 in \cite{HS01}).

\medskip

Thus, the spectrum of the magnetic Schr\"odinger operator $H_{\a}$ on the periodic graph $\cG$ is given by
\[\lb{spec}
\s(H_{\a})=\bigcup_{\vt\in\T^d}\s\big(H_{\gm,\f}(\vt)\big)=
\bigcup_{n=1}^{\nu}\s_n(H_{\a}).
\]
Note that if $\l_{\a,n}(\cdot)= C_{n}=\const$ on some subset of $\T^d$ of
positive Lebesgue measure, then  the operator $H_\a$ on $\cG$ has the
eigenvalue $C_{n}$ of infinite multiplicity. We call $C_{n}$
a \emph{flat band}.

Thus, the spectrum of the magnetic Schr\"odinger operator
$H_{\a}$ on the periodic graph $\cG$ has the form
$$
\s(H_{\a})=\s_{ac}(H_{\a})\cup \s_{fb}(H_{\a}),
$$
where $\s_{ac}(H_{\a})$ is the absolutely continuous spectrum, which is a
union of non-degenerate bands, and $\s_{fb}(H_{\a})$ is the set of
all flat bands (eigenvalues of infinite multiplicity). An open
interval between two neighboring non-degenerate bands is
called a \emph{spectral gap}.

\subsection{Estimates of the Lebesgue measure of the spectrum}
Now we estimate the position of the bands $\s_n(H_{\a})$, $n\in\N_\n$,
defined by \er{ban.1} in terms of eigenvalues of the magnetic Schr\"odinger operator on some subgraph of the fundamental graph $G_*$ and the maximal vertex degree of the remaining subgraph of $G_*$. We also obtain the estimate of the Lebesgue measure of the spectrum of the magnetic Schr\"odinger operator $H_{\a}$ on a periodic graph $\cG$ in terms of the invariant $\cI$ given in \er{dIm}.

Let $\gm\in\mF(\k)$ be a minimal form on the fundamental graph $G_*=(\cV_*,\cE_*)$, where $\mF(\k)$ is given by \er{vv1f} at $\mathbf{x}=\k$, and $\k$ is the coordinate form defined by \er{edco}, \er{dco}. In order to formulate our result we equip each unoriented edge $\be\in\cE_*$ of $G_*$ with some orientation and represent $G_*$ as a union of two graphs with the same vertex set $\cV_*$:
\[\lb{cEgm}
\begin{aligned}
G_*=G_\gm^0\cup G_\gm, \qqq G_\gm=(\cV_*,\cE_\gm),\qqq G_\gm^0=(\cV_*,\cE_*\sm\cE_\gm),\\ \textrm{where}\qqq \cE_\gm=\cE_*\cap\supp \gm.
\end{aligned}
\]
The graph $G_\gm^0$ is obtained from the fundamental graph $G_*$ by deleting all edges of $\supp \gm$ and preserving the vertex set $\cV_*$. Let $\D_{\a}^0$ be the magnetic Laplacian defined by \er{DLO} on the graph $G_\gm^0$. Since the graph $G_\gm^0$ consists of $\n=\#\cV_*$ vertices, the magnetic Schr\"odinger operator $H_{\a}^0=\D_{\a}^0+Q$ on the graph $G_\gm^0$ has $\n$ eigenvalues $\m_{\a,n}$, $n\in\N_\n$, which are labeled (counting multiplicities) by
\[\lb{evH0}
\m_{\a,1}\leq\ldots\leq\m_{\a,\n}.
\]

\begin{theorem}
\lb{T1} i) Each band $\s_n(H_{\a})$, $n\in\N_\n$, of the magnetic Schr\"odinger operator $H_{\a}=\D_{\a}+Q$ on a periodic graph $\cG$
defined by \er{ban.1} satisfies
\[\lb{fesbp1}
\s_n(H_{\a})\ss\big[\m_{\a,n},\m_{\a,n}+2\vk^{\gm}_+\big],
\qq \textrm{where}\qq
\vk^{\gm}_+=\max\limits_{v\in \cV_*}\vk_v^{\gm},
\]
$\gm\in\mF(\k)$ is a minimal form, $\mu_{\a,n}$ are given by (\ref{evH0}), and  $\vk_v^{\gm}$ is the degree of the vertex $v\in\cV_*$ on the graph $G_\gm$ defined in \er{cEgm}.

ii) The Lebesgue measure $|\s(H_{\a})|$ of the spectrum of $H_{\a}$ satisfies
\[
\lb{eq.7'}
|\s(H_{\a})|\le \sum_{n=1}^{\n}|\s_n(H_{\a})|\le 4\cI,
\]
where the invariant $\cI$ is given in \er{dIm}.

iii) There exists a periodic graph $\cG$ such that the estimate \er{eq.7'} becomes an identity for any potentials $Q$ and $\a$.
\end{theorem}

\no \textbf{Remarks.} 1) The decomposition of the magnetic Schr\"odinger operator $H_{\a}$ into the direct integral \er{raz1m} -- \er{l2.13am} in terms of minimal forms is essentially used in the proof of the estimate \er{eq.7'}.

2) From Theorem \ref{T1}.\emph{iii}) it follows that
the Lebesgue measure of the spectrum of $H_{\a}$ can be arbitrarily large (for specific graphs).

3) In \cite{KS17} the authors also estimated the Lebesgue measure of the spectrum of the magnetic Schr\"odinger operators $H_{\a}$ with periodic magnetic and electric potentials on periodic graphs in terms of the Betti number $\b$:
\[
\lb{esLm1}
|\s(H_{\a})|\le \sum_{n=1}^{\n}|\s_n(H_{\a})|\le 4\b.
\]
Due to Proposition \ref{TNNI0}.\emph{iii}), the difference between the Betti number $\b$ and the invariant $\cI$ given in \er{dIm} can be arbitrarily large. Then comparing estimates \er{eq.7'} and  \er{esLm1} we see that the estimate \er{eq.7'} is better than \er{esLm1}.

\medskip

Now we estimate a variation of the spectrum of the Schr\"odinger operators under a perturbation by magnetic fields.

\begin{theorem}
\lb{Temf} Consider the magnetic Schr\"odinger operator $H_\a=\D_\a+Q$ and the Schr\"odinger operator $H_0=\D_0+Q$ without a magnetic field. Then the spectral bands
$$
\s_{\a,n}=\s_{n}(H_{\a})=[\l_{\a,n}^-,\l_{\a,n}^+]\qqq\textrm{and}\qqq \s_{0,n}=\s_{n}(H_0)=[\l_{0,n}^-,\l_{0,n}^+]
$$
and their band ends $\l_{\a,n}^\pm$ and $\l_{0,n}^\pm$ satisfy
\[
\lb{emf1} \L_1\le \l_{\a,n}^\pm-\l_{0,n}^\pm\le \L_\n,
\]
\[
\lb{emf2} \big||\s_{\a,n}|-|\s_{0,n}|\big|\le\L_\n-\L_1,
\]
where
\[\lb{emf3'}
\begin{array}{c}
\L_1=\min\limits_{\vt\in\T^d} \l_1(X_{\gm,\wt\f}(\vt)),\qqq
\L_\n=\max\limits_{\vt\in\T^d} \l_\n(X_{\gm,\wt\f}(\vt)),\\[12pt] X_{\gm,\wt\f}(\cdot)=\D_{\gm,\wt\f}(\cdot)-\D_{\gm,0}(\cdot),
\end{array}
\]
$\D_{\gm,\wt\f}(\cdot)$ is given by \er{mmfr}, and $\wt\phi$ is defined by \er{wtcat}. Moreover,
$\L_1$ and $\L_\n$ satisfy the following estimates:
\[
\begin{aligned}
\lb{es30} \max\{|\L_1|,|\L_\n|\}\le C_{\wt\f},\qqq \L_\n-\L_1\le
2C_{\wt\f},
\end{aligned}
\]
where
\[
\begin{aligned}
\lb{es3} C_{\wt\f}=2\max_{v\in \cV_*}\sum_{\be=(v,u)\in\supp\wt\phi}\big|\sin\textstyle{\wt\phi(\be)\/2}\big|.
\end{aligned}
\]
\end{theorem}

\no {\bf Remark.} In \cite{KS17} the authors estimated a variation of the spectrum of the Schr\"odinger operators under a perturbation by a magnetic field in terms of magnetic fluxes:
\[\lb{emf22}
\big||\s_{\a,n}|-|\s_{0,n}|\big|\le2C_{\a_*}, \qqq\textrm{where}\qqq
C_{\a_*}=2\max_{v\in \cV_*}\sum_{\be=(v,u)\in\supp\a_*}\big|\sin\textstyle{\a_*(\be)\/2}\big|,
\]
where $\a_*\in\mF(\a)$ is a 1-form (in general, not minimal) equal to zero on edges of a spanning tree of the fundamental graph $G_*$ and coinciding with fluxes of $\a$ on other edges of $G_*$. Due to the definition \er{wtcat} of the form $\wt\phi$ and Remarks 1,2) on page \pageref{Rem2}, in the general case, we have $\#\supp\wt\phi<\#\supp\a_*$. Moreover, in the case $\cI_{\gm,\f}=d$, where $\cI_{\gm,\f}$ is defined by \er{decI}, the identity \er{es3} immediately gives that $C_{\wt\f}=0$ and, consequently, $\s(H_\a)=\s(H_0)$. This result cannot be deduced from the estimate \er{emf22}. Thus, in general, the estimate \er{emf2} is better than \er{emf22}.

\medskip

The paper is organized as follows. In Section \ref{Sec2} we give a full description of minimal forms on a finite connected graph (see Theorem \ref{Pphi}).

Section \ref{Sec3} is devoted to the direct integral decomposition for the magnetic Schr\"odinger operators $H_\a$ with magnetic vector potentials $\a$ on periodic graphs. In this section we consider some properties of the coordinate form $\k$ on the fundamental graph $G_*$ defined by \er{edco}, \er{dco} (see Proposition \ref{Pal0}) and give a general representation of fiber magnetic Laplacians in terms of any 1-forms $(\gb,\ga)\in\mF(\k)\ts\mF(\a)$ (see Theorem \ref{TDIBB}).  In Section \ref{Sec3} we also prove Theorem \ref{TDImf} about the direct integral for the Schr\"odinger operators, where fiber Laplacians have the minimal numbers of coefficients depending on the quasimomentum and on the magnetic potential, Proposition \ref{TNNI0} about some properties of the invariants $\cI$ and $\cI_\a$ given by \er{dIm} and Corollary \ref{TCo0}.

In Section \ref{Sec4}, using the representation of fiber magnetic Laplacians in terms of minimal forms, we prove Theorem \ref{T1} about spectral estimates for the magnetic Schr\"odinger operators and show that these estimates become identities for specific graphs (see Proposition \ref{TG1}). In this section we also prove Theorem \ref{Temf} about the estimates of a variation of the spectrum of the Schr\"odinger operators under a perturbation by a magnetic field.

\section{\lb{Sec2} Minimal forms and their properties}
\setcounter{equation}{0}

The results of this section repeat the results of Section 3 in \cite{KS18}. For the convenience of the reader, we include them also here.

\subsection{Betti numbers and spanning trees}
We recall the definitions of the Betti number and spanning trees which will be used in the proof of our results. Let $G=(V,E)$ be a finite connected graph. Recall that $\#\cM$ denotes the number of elements in a set $\cM$.

\medskip

$\bullet$  A \emph{spanning tree} $T=(V,E_{T})$ of the graph $G$ is a connected subgraph of $G$ which has no cycles and contains all vertices of $G$.

$\bullet$ The \emph{Betti number} $\b$ of the graph $G$ is defined as
\[\lb{benu}
\b=\# E-\# V+1.
\]
Note that the Betti number $\b$ can also be defined in one of the following equivalent ways:
\begin{itemize}
  \item[\emph{i})] as the number of edges that have to be removed from $E$ to turn $G$ into a spanning tree of $G$;
  \item[\emph{ii})] as the dimension of the cycle space $\cC$ of the graph $G$, i.e.,
\[
\lb{debn}
\b=\dim\cC.
\]
\end{itemize}

\medskip

We introduce the set $\cS_T$ of all edges from $E$ that do not belong to the spanning tree $T$, i.e.,
\[
\lb{ulS}
\cS_T=E\sm E_T,
\]
and recall some properties of spanning trees (see, e.g., Lemma 5.1 and Theorem 5.2 in \cite{B74}, and Proposition 1.3 in \cite{CDS95}).

\medskip

\no\textbf{Properties of spanning trees.}\lb{PSTs}

\emph{1) The set $\cS_T$ consists of exactly $\b$ edges, where $\b$ is the Betti number defined by \er{benu}.}

\emph{2) For each $\be\in\cS_T$ there exists a unique cycle $\mathbf{c}_\be$
consisting of only $\be$ and edges of $T$.}

\emph{3) The set of all such cycles $\mathbf{c}_\be$, $\be\in\cS_T$, forms a basis $\cB_T$ of the cycle space $\cC$ of the graph $G$:}
\[\lb{bst}
\cB_T=\{\mathbf{c}_\be: \be\in\cS_T\}.
\]

\emph{4) Let $0=\x_1<\x_2\leq\x_3\leq\ldots\leq\x_\n$ be the eigenvalues of the Laplacian $\D$ on the graph $G=(V,E)$, where $\n=\# V$. Then the number of spanning trees of $G$ is equal to $\frac1\n\,\x_2\,\x_3\ldots\x_\n$.}

\medskip

\no \textbf{Example.} For the graph $G$ shown in Fig.\ref{ffS'}\emph{a} we can choose the spanning trees $T$ and $\wt T$ (Fig.\ref{ffS'} \emph{b,c}).
The set $\cS_T$ consists of three edges $\be_1,\be_2,\be_3$ (they are shown in Fig.\ref{ffS'} \emph{b,c} by the dotted lines) and depends on the choice of the spanning tree. The Betti number $\b$ defined by \er{benu} is equal to 3 and does not depend on the set $\cS_T$.

\setlength{\unitlength}{1.0mm}
\begin{figure}[h]
\centering
\unitlength 1.0mm % = 2.845pt
\linethickness{0.4pt}
%\ifx\plotpoint\undefined\newsavebox{\plotpoint}\fi % GNUPLOT compatibility
\begin{picture}(100,30)

\put(10,10){\circle{1}}
\put(10,30){\circle{1}}
\put(10,20){\circle{1}}
\put(0,20){\circle{1}}
\put(20,20){\circle{1}}
\put(10,10){\line(1,1){10.00}}
\put(10,10){\line(-1,1){10.00}}
\put(10,30){\line(1,-1){10.00}}
\put(10,30){\line(-1,-1){10.00}}

\put(10,30){\line(0,-1){10.00}}
\put(0,20){\line(1,0){20.00}}

%**************************************
\put(50,10){\circle{1}}
\put(50,30){\circle{1}}
\put(50,20){\circle{1}}
\put(40,20){\circle{1}}
\put(60,20){\circle{1}}
\put(50,10){\line(1,1){10.00}}
\put(50,10.2){\line(1,1){10.00}}
\put(50,9.8){\line(1,1){10.00}}
\put(50,10.1){\line(1,1){10.00}}
\put(50,9.9){\line(1,1){10.00}}

\put(50,30){\line(0,-1){10.00}}
\put(50.1,30){\line(0,-1){10.00}}
\put(49.9,30){\line(0,-1){10.00}}
\put(40,20.1){\line(1,0){20.00}}
\put(40,19.9){\line(1,0){20.00}}
\put(40,20){\line(1,0){20.00}}

\qbezier[20](50,10)(45,15)(40,20)
\qbezier[20](50,30)(55,25)(60,20)
\qbezier[20](50,30)(45,25)(40,20)

%\put(50,10.2){\line(-1,1){10}}
%\put(50,9.8){\line(-1,1){10}}
%\multiput(50,30)(4,-4){3}{\line(1,-1){2}}
%\multiput(50,30)(-4,-4){3}{\line(-1,-1){2}}
\put(42,13){$\be_1$}
\put(41,26){$\be_2$}
\put(55,26){$\be_3$}
%***********************************
\put(90,10){\circle{1}}
\put(90,30){\circle{1}}
\put(90,20){\circle{1}}
\put(80,20){\circle{1}}
\put(100,20){\circle{1}}
\put(90,10){\line(1,1){10.00}}
\put(90,10.2){\line(1,1){10.00}}
\put(90,9.8){\line(1,1){10.00}}
\put(90,10.1){\line(1,1){10.00}}
\put(90,9.9){\line(1,1){10.00}}

\put(90,30){\line(0,-1){10.00}}
\put(90.1,30){\line(0,-1){10.00}}
\put(89.9,30){\line(0,-1){10.00}}
\put(90,20.1){\line(1,0){10.00}}
\put(90,19.9){\line(1,0){10.00}}
\put(90,20){\line(1,0){10.00}}

\qbezier[15](80,20)(85,20)(90,20)
\qbezier[20](90,10)(85,15)(80,20)
\qbezier[20](90,30)(95,25)(100,20)

%\multiput(80,20)(4,0){3}{\line(1,0){2}}
%\multiput(90,10)(-4,4){3}{\line(-1,1){2}}
%\multiput(90,30)(4,-4){3}{\line(1,-1){2}}
\put(90,30.2){\line(-1,-1){10.00}}
\put(90,29.8){\line(-1,-1){10.00}}
\put(90,30.1){\line(-1,-1){10.00}}
\put(90,29.9){\line(-1,-1){10.00}}
\put(90,30){\line(-1,-1){10.00}}
%\multiput(90,30)(-4,-4){3}{\line(-1,-1){2}}
\put(82,13){$\be_1$}
\put(84,21){$\be_2$}
\put(95,26){$\be_3$}

\put(-3,10){(\emph{a})}
%**********************************
\put(33,10){(\emph{b})}
\put(73,10){(\emph{c})}

\put(0,26){$G$}
\put(33,26){$T$}
\put(80,26){$\wt T$}
\end{picture}

\vspace{-10mm}
\caption{\footnotesize  \emph{a}) A finite connected graph $G$;\quad \emph{b}),\,\emph{c}) the spanning trees $T$ and $\wt T$, $\cS_T=\{\be_1,\be_2,\be_3\}$, $\b=3$.} \label{ffS'}
\end{figure}
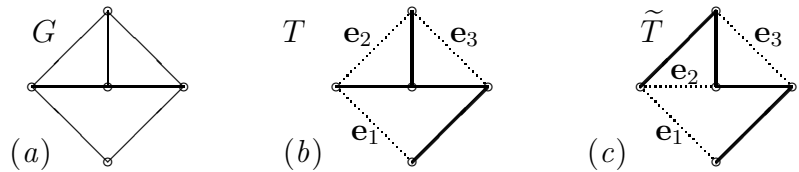

\subsection{Minimal forms}
Let $\mathbf{x}:A\to\R^n$, $n\in\N$, be a 1-form on a finite connected graph $G=(V,E)$, where $A$ is the set of oriented edges of $G$. We describe all minimal forms $\gm\in\mF(\mathbf{x})$ on $G$, where $\mF(\mathbf{x})$ is defined by \er{vv1f}, and give an explicit expression for the number of edges in their supports.

Let $T=(V,E_{T})$ be a spanning tree of the graph $G$. Then, due to the properties 2) and 3) of spanning trees (see page \pageref{PSTs}), for each $\be\in\cS_T=E\sm E_T$, there exists a unique cycle $\mathbf{c}_\be$ consisting of only $\be$ and edges of $T$ and the set of all these cycles $\cB_T$ forms a basis of the cycle space $\cC$ of the graph $G$.

Among all spanning trees of the graph $G$ we specify a subset of \emph{minimal} spanning trees. Let $\b_T=\b_T(\mathbf{x})$ be the number of basic cycles from $\cB_T$ with non-zero fluxes of the form $\mathbf{x}$:
\[\lb{ncIT}
\b_T=\#\{\mathbf{c}\in\cB_T : \Phi_{\mathbf{x}}(\mathbf{c})\neq0\}.
\]
A spanning tree $T$ of the graph $G$ is called \emph{$\mathbf{x}$-minimal} if
\[\lb{dmf}
\b_T\leq\b_{\wt T} \qqq \textrm{for any spanning tree $\wt T$ of $G$.}
\]
Among all spanning trees $T$ of $G$ the $\mathbf{x}$-minimal one gives a minimal number of basic cycles from $\cB_T$ with non-zero fluxes of the form $\mathbf{x}$. Since the set $\cB_T$ of basic cycles is finite, minimal spanning trees exist. All $\mathbf{x}$-minimal spanning trees $T$ have the same number $\b_T$. We denote this number by $\b(\mathbf{x})$:
\[\lb{nbo}
\b(\mathbf{x})=\b_T \qqq \textrm{for any $\mathbf{x}$-minimal spanning tree $T$ of $G$.}
\]

\medskip

\no \textbf{Example.} We consider the graph $G$ shown in Fig.\ref{FkfKl}\emph{a}. The values of the 1-form $\mathbf{x}:A\to\R^2$ on $G$ are shown in the figure near the corresponding edges. For the spanning tree $T$ shown in Fig.\ref{FkfKl}\emph{b} the set $\cS_T$ consists of 4 edges $\be_2,\be_4,\be_5,\be_6$ (they are shown by the dotted lines) and the set $\cB_T$ defined by \er{bst} consists of 4 basic cycles. Only one of them, the cycle $(\be_1,\be_2,\be_3)$, has zero flux of the form $\mathbf{x}$. Thus, the number $\b_T$ defined by \er{ncIT} is equal to $3$. For the spanning tree $\wt T$ (Fig.\ref{FkfKl}\emph{c}) all basic cycles from $\cB_{\wt T}$ have non-zero fluxes of the form $\mathbf{x}$. Thus, $\b_{\wt T}=4$. Looking over all spanning trees of $G$ (the number of spanning trees is finite) one can check that $T$ is a $\mathbf{x}$-minimal spanning tree of the graph $G$.

\begin{figure}[h]
\centering
\unitlength 0.9mm % = 2.845pt
\linethickness{0.4pt}
\ifx\plotpoint\undefined\newsavebox{\plotpoint}\fi % GNUPLOT compatibility
\hspace{-2.5cm}
\begin{picture}(146,48)(0,0)

%***************************
\put(20,6){(\emph{a})}
\put(25,38){$G$}
\put(31.0,27.0){\vector(1,2){1.0}}
\put(44.2,28.5){\vector(1,-2){1.0}}
\put(38,15.0){\vector(-1,0){1.0}}

\put(25,15){\line(1,0){26.0}}
\put(25,15){\line(1,2){13.0}}

\put(50.8,29.5){\vector(-1,2){1.0}}
\put(26.2,31.5){\vector(-1,-2){1.0}}
\put(38,9.0){\vector(1,0){1.0}}

\put(51,15){\line(-1,2){13.0}}

\put(25,15){\circle*{1.5}}
%\put(19.0,13.5){$v_1$}
\put(51,15){\circle*{1.5}}
%\put(52.5,13.5){$v_3$}
\put(38,41){\circle*{1.5}}
%\put(38.0,42){$v_2$}
\bezier{200}(25,15)(38,3)(51,15)
\bezier{200}(25,15)(20,33)(38,41)
\bezier{200}(51,15)(56,33)(38,41)

\put(32.0,17.0){$\scriptstyle(-{1\/2},0)$}
\put(32.0,27.0){$\scriptstyle(0,{1\/2})$}
\put(36.5,22.0){$\scriptstyle({1\/2},-{1\/2})$}

\put(36.0,5.0){$\scriptstyle(-{1\/2},0)$}
\put(17.0,30.0){$\scriptstyle(0,{1\/2})$}
\put(51.0,30.0){$\scriptstyle({1\/2},-{1\/2})$}

%***************************************
%\put(95,46){the minimal forms $\gm$}
\put(70,6){(\emph{b})}
\put(75,38){$T$}
%\put(75,15){\line(1,0){26.0}}
\bezier{40}(75,15)(88,15)(101,15)
\put(88,15.0){\vector(-1,0){1.0}}

%\put(75,15){\line(1,2){13.0}}
\bezier{40}(75,15)(81.5,28)(88,41)
\put(81.2,27.2){\vector(1,2){1.0}}

\put(94.2,28.5){\vector(1,-2){1.0}}

\put(100.8,29.5){\vector(-1,2){1.0}}
\put(100.5,29.5){\vector(-1,2){1.0}}

\put(76.2,31.5){\vector(-1,-2){1.0}}
\put(76.5,31.5){\vector(-1,-2){1.0}}

%\put(88,9.3){\vector(1,0){1.0}}
\put(88,9.0){\vector(1,0){1.0}}
%\put(88,8.8){\vector(1,0){1.0}}

%\put(101,15){\line(-1,2){13.0}}
\bezier{40}(101,15)(94.5,28)(88,41)

%\bezier{200}(75,15.4)(88,3.4)(101,15.4)
%\bezier{200}(75,15.2)(88,3.2)(101,15.2)
\bezier{50}(75,15)(88,3)(101,15)
%\bezier{200}(74.8,15)(87.8,3)(100.8,15)
%\bezier{200}(74.6,15)(87.6,3)(100.6,15)

\bezier{200}(75,15)(70,33)(88,41)
\bezier{200}(75.2,15)(70.2,33)(88.2,40.8)
%\bezier{200}(75.4,15)(70.4,33)(88.4,40.7)
%\bezier{200}(74.8,15)(69.8,33)(87.8,41.2)

\bezier{200}(101,15)(106,33)(88,41)
\bezier{200}(100.8,15)(105.8,33)(87.8,40.7)
%\bezier{200}(100.6,15)(105.6,33)(87.6,41)
%\bezier{200}(101.2,15)(106.2,33)(88.2,41.2)

\put(75,15){\circle*{1.5}}
%\put(69.0,13.5){$v_1$}
\put(101,15){\circle*{1.5}}
%\put(102.5,13.5){$v_3$}
\put(88,41){\circle*{1.5}}
%\put(88.0,42){$v_2$}

\put(86.0,17.0){$\be_5$}
\put(81.0,24.0){$\be_6$}
\put(91.0,24.0){$\be_4$}
\put(86.0,5.0){$\be_2$}
\put(71.0,32.0){$\be_1$}
\put(100.0,32.0){$\be_3$}

%***************************************
\put(147,38){$\wt T$}
\put(120,6){(\emph{c})}
%\put(125,14.6){\line(1,0){26.0}}
\put(125,14.8){\line(1,0){26.0}}
\put(125,15){\line(1,0){26.0}}
\put(125,15.2){\line(1,0){26.0}}
%\put(125,15.4){\line(1,0){26.0}}

\put(138,15.0){\vector(-1,0){1.0}}
\put(138,14.8){\vector(-1,0){1.0}}
\put(138,15.2){\vector(-1,0){1.0}}

\bezier{40}(125,15)(131.5,28)(138,41)

\put(131.2,27.5){\vector(1,2){1.0}}

%\put(144.4,28.5){\vector(1,-2){1.0}}
\put(144.2,28.5){\vector(1,-2){1.0}}
%\put(144.0,28.5){\vector(1,-2){1.0}}

\put(150.8,29.5){\vector(-1,2){1.0}}
\put(150.6,29.5){\vector(-1,2){1.0}}

\put(126.2,31.5){\vector(-1,-2){1.0}}
%\put(126.5,31.5){\vector(-1,-2){1.0}}

\put(138,9.0){\vector(1,0){1.0}}

%\put(151,14.2){\line(-1,2){13.0}}
%\put(151,14.6){\line(-1,2){13.0}}
%\put(151,15){\line(-1,2){13.0}}
%\put(151,15.4){\line(-1,2){13.0}}
%\put(151,15.8){\line(-1,2){13.0}}
\bezier{40}(151,15)(144.5,28)(138,41)

\bezier{50}(125,15)(138,3)(151,15)
\bezier{50}(125,15)(120,33)(138,41)

\bezier{200}(151,15)(156,33)(138,41)
\bezier{200}(150.8,15)(155.8,33)(137.8,40.7)
%\bezier{200}(150.6,15)(155.6,33)(137.6,41)
%\bezier{200}(151.2,15)(156.2,33)(138.2,41.2)

\put(125,15){\circle*{1.5}}
%\put(119.0,13.5){$v_1$}
\put(151,15){\circle*{1.5}}
%\put(152.5,13.5){$v_3$}
\put(138,41){\circle*{1.5}}
%\put(138.0,42){$v_2$}

\put(136.0,17.0){$\be_5$}
\put(131.0,24.0){$\be_6$}
\put(141.0,24.0){$\be_4$}
\put(136.0,5.0){$\be_2$}
\put(121.0,32.0){$\be_1$}
\put(150.0,32.0){$\be_3$}

%\put(135.0,16.5){$\scriptstyle(-1,0)$}
%\put(132.0,27.0){$\scriptstyle(0,1)$}
%\put(136.5,22.5){$\scriptstyle(1,-1)$}
%
%\put(136.0,5.0){$\scriptstyle(0,0)$}
%\put(118.0,32.0){$\scriptstyle(0,0)$}
%\put(150.0,32.0){$\scriptstyle(0,0)$}
%\put(70,-2){$d=2$, \qq $\cI=3$, \qq $\b=4$}
\end{picture}

%\vspace{-0.7cm}
\caption{\footnotesize  \emph{a}) A graph $G$ with the 1-form $\mathbf{x}$;\quad \emph{b}) a $\mathbf{x}$-minimal spanning tree $T$;\quad \emph{c}) a non-minimal spanning tree $\wt T$.}
\label{FkfKl}
\end{figure}
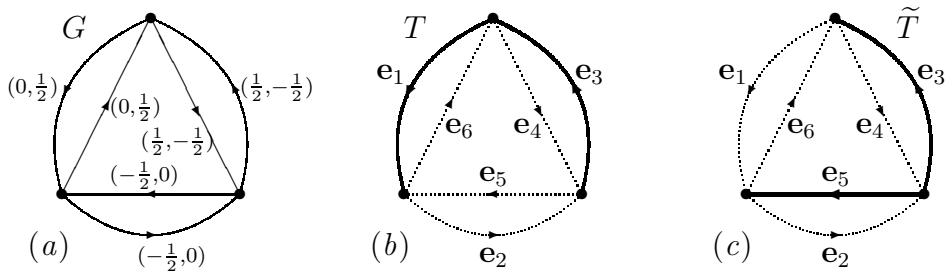

Let $T=(V,E_T)$ be a spanning tree of the graph $G$. We equip each edge of $\cS_T=E\sm E_T$ with some orientation and associate with $T$ a 1-form $\gm(\,\cdot\,,T): A\to\R^n$ equal to zero on edges of $T$ and coinciding with the flux $\Phi_{\mathbf{x}}(\mathbf{c}_\be)$ of the form $\mathbf{x}$ through the cycle $\mathbf{c}_\be$ on each edge $\be\in\cS_T$, i.e.,
\[\lb{cat}
\gm(\be,T)=\left\{
\begin{array}{cl}
 \Phi_{\mathbf{x}}(\mathbf{c}_\be), &  \textrm{ if } \, \be\in\cS_T\\[6pt]
 -\Phi_{\mathbf{x}}(\mathbf{c}_{\ul\be}), &  \textrm{ if } \, \ul\be\in\cS_T\\[6pt]
  0, \qq & \textrm{ otherwise} \\
\end{array}\right..
\]

\begin{theorem}
\lb{Pphi} Let $\mathbf{x}:A\to\R^n$, $n\in\N$, be a 1-form on a finite connected graph $G=(V,E)$, and $\mF(\mathbf{x})$ be given by \er{vv1f}. Then

i) Each minimal form $\wt\gm\in\mF(\mathbf{x})$ satisfies $\wt\gm=\gm(\,\cdot\,,T)$, where the function $\gm(\,\cdot\,,T): A\ra\R^n$ is defined by \er{cat} for some $\mathbf{x}$-minimal spanning tree $T$ of the graph $G$.

ii) Conversely, for each $\mathbf{x}$-minimal spanning tree $T$ of the graph $G$ the function $\gm(\,\cdot\,,T)$ is a minimal form in the set $\mF(\mathbf{x})$.

iii) The number $\b(\mathbf{x})$ defined by \er{ncIT}, \er{nbo}
satisfies
\[\lb{gfe1}
\b(\mathbf{x})=\textstyle\frac12\,\#\supp \gm,
\]
where $\gm\in\mF(\mathbf{x})$ is a minimal form.

iv) There exist a finite connected graph $\wt G$ and a 1-form $\wt{\mathbf{x}}$ on it such that $\gm(\,\cdot\,,T)=\gm(\,\cdot\,,\wt T\,)$ for some distinct $\wt{\mathbf{x}}$-minimal spanning trees $T$ and $\wt T$ of the graph $\wt G$.

v) Let the function $\gm(\,\cdot\,,T)$ be defined by \er{cat} for some spanning tree $T$ of the graph $G$. If $\gm(\be,T)$, $\be\in\supp\gm(\,\cdot\,,T)\cap\cS_T$, are rationally independent, then $\gm(\,\cdot\,,T)$ is a minimal form.

\end{theorem}

\no{\bf Proof.} Items \emph{i}) -- \emph{iv)} were proved in \cite{KS18} (Theorem 3.1).

\emph{v)} The proof is by contradiction. Let $\wt\gm\in\mF(\mathbf{x})$ be a minimal form and
\[\lb{wtgm}
\#\supp\wt\gm<\#\supp\gm(\,\cdot\,,T).
\]
Then, due to item \emph{i}), there exists a $\mathbf{x}$-minimal spanning tree $\wt T$ of the graph $G$ such that
\[\lb{gmgb1}
\wt\gm=\gm(\,\cdot\,,\wt T).
\]
Combining \er{wtgm} and \er{gmgb1}, we obtain
\[\lb{wtTT}
\#\supp\gm(\,\cdot\,,\wt T)<\#\supp\gm(\,\cdot\,,T).
\]
Due to the properties 2) and 3) of spanning trees (see page \pageref{PSTs}), for each $\be\in\cS_T=E\sm E_T$, there exists a unique cycle $\mathbf{c}_\be$ consisting of only $\be$ and edges of $T$ and the set $\cB_T$ of all these cycles forms a basis of the cycle space $\cC$ of the graph $G$. Due to the definition \er{cat} of the form $\gm(\,\cdot\,,T)$,
$$
\textstyle\frac12\#\supp\gm(\,\cdot\,,T)=\b_T, \qqq \b_T=\#\{\mathbf{c}\in\cB_T : \Phi_{\mathbf{x}}(\mathbf{c})\neq0\}.
$$
This and \er{wtTT} yield $\b_{\wt T}<\b_T$. Consequently, there exists a cycle $\mathbf{c}\in\cB_{\wt T}\cap\ker\F_\mathbf{x}$ satisfying the following conditions

$\bu$ $\mathbf{c}=\sum\limits_{\be\in\cS_T}a_\be\mathbf{c}_\be$, for some $a_\be\in\Z$,

$\bu$ there exists $\be\in\cS_T$ such that $\F_\mathbf{x}(\mathbf{c}_\be)\neq0$ and $a_\be\neq0$.\\
From this and the identity $\gm(\be,T)=\F_\mathbf{x}(\mathbf{c}_\be)$, $\be\in\cS_T$, we deduce that
$$
0=\F_\mathbf{x}(\mathbf{c})=\sum\limits_{\be\in\cS_T \atop
\mathbf{c}_\be\notin\ker\F_\mathbf{x}}a_\be\F_\mathbf{x}(\mathbf{c}_\be), \qq \textrm{where $\F_\mathbf{x}(\mathbf{c}_\be)$ are rationally independent},
$$
which is a contradiction. Thus, $\gm(\,\cdot\,,T)$ is a minimal form. \qq \BBox

\medskip

\no {\bf Remarks.} 1) Theorem \ref{Pphi}.\emph{i}) -- \emph{ii}) gives a full description of minimal forms $\gm\in\mF(\mathbf{x})$. Each minimal form $\gm$ is given by \er{cat} for some $\mathbf{x}$-minimal spanning tree $T$ and, conversely, for each $\mathbf{x}$-minimal spanning tree $T$ the form $\gm(\,\cdot\,,T)$ defined by \er{cat} is minimal.

2)  Due to Theorem \ref{Pphi}.\emph{iv}), the correspondence between minimal forms in a fixed set $\mF(\mathbf{x})$ on a graph $G$ and $\mathbf{x}$-minimal spanning trees of $G$ is not a bijection. Two different spanning trees may be associated with the same minimal form. Thus, the number of minimal forms in the set $\mF(\mathbf{x})$ is not greater than the number of $\mathbf{x}$-minimal spanning trees of the graph $G$.

\section{\lb{Sec3} Direct integrals}
\setcounter{equation}{0}

In this section we consider some properties of the coordinate form $\k$ on the fundamental graph $G_*$ defined by \er{edco}, \er{dco} and introduce an important 1-form $\t\in\mF(\k)$ on the graph $G_*$ called the \emph{index form}. We also give a general representation of fiber operators for the magnetic Laplacians $\D_\a$ with magnetic vector potentials $\a$ in terms of any 1-forms $(\gb,\ga)\in\mF(\k)\ts\mF(\a)$. Then we prove Theorem \ref{TDImf} about the decomposition for the magnetic Schr\"odinger operators in terms of minimal forms and Proposition \ref{TNNI0} about some properties of the invariants $\cI$ and $\cI_\a$ given by \er{dIm}.

\subsection{Properties of coordinate forms}
We formulate some properties of the coordinate form $\k$ defined by \er{edco}, \er{dco} on the fundamental graph $G_*$. Recall that $\cC$ is the cycle space of $G_*$.

\begin{proposition}
\lb{Pal0}
i) The image of the flux function $\Phi_{\k}:\cC\to\R^d$ of the coordinate form $\k:\cA_*\ra\R^d$ defined by \er{mafla} is the lattice $\Z^d$:
\[\lb{prao}
\Phi_{\k}(\cC)=\Z^d.
\]

ii) The kernel of the flux function $\Phi_{\k}$ does not depend on the choice of

$\bu$ the embedding of $\cG$ into $\R^d$;

$\bu$ the basis $a_1,\ldots,a_d$ of the lattice $\G$,

\no i.e., $\ker\Phi_{\k}$
is an invariant of the periodic graph $\cG$.

iii) The dimension of the kernel of $\Phi_{\k}$ satisfies
\[\lb{dike}
\dim\ker\Phi_{\k}=\b-d,
\]
where $\b$ is the Betti number of the fundamental graph $G_*$ given in \er{nni0}.
\end{proposition}

This proposition was proved in \cite{KS18} (see Proposition 4.1).

\medskip

\no \textbf{Remark.} Any cycle $\textbf{c}$ on the fundamental graph $G_*$ is obtained by factorization of a path on the periodic graph $\cG$ connecting some $\G$-equivalent vertices $v\in\cV$ and $v+a\in\cV$, $a\in\G$. Furthermore, the flux $\Phi_{\k}(\textbf{c})$ of the coordinate form $\k$ through the cycle $\mathbf{c}$ is equal to $\mn=(n_1,\ldots,n_d)\in\Z^d$, where $a=n_1a_1+\ldots+n_da_d$. In particular, $\Phi_{\k}(\textbf{c})=0$ if and only if the cycle $\textbf{c}$ on $G_*$ corresponds to a cycle on $\cG$.

\subsection{Edge indices.} In order to prove Theorem \ref{TDImf} about the direct integral decomposition we need to define an {\it edge index}, which was introduced in \cite{KS14}. The indices are important to study the spectrum of the Laplacians and Schr\"odinger operators on periodic graphs, since the initial fiber operators are expressed in terms of indices of the fundamental graph edges (see \er{l2.15''}).

For any vertex $v\in\cV$ of a $\G$-periodic graph $\cG$ the following unique representation holds true:
\[
\lb{Dv} v=\{v\}+[v], \qq \textrm{where}\qq \{v\}\in \Omega,\qquad [v]\in\G,
\]
$\Omega$ is the fundamental cell of the lattice $\G$ defined by \er{fuce}.
In other words, each vertex $v$ can be obtained from a vertex $\{v\}\in \Omega$ by a shift by the vector $[v]\in\G$. We call $\{v\}$ and $[v]$ the \emph{fractional and integer parts of the vertex $v$}, respectively.
For any oriented edge $\be=(u,v)\in\cA$ we define the {\bf edge "index"}
$\t(\be)$ as the vector of the lattice $\Z^d$ given by
\[
\lb{in}
\t(\be)=[v]_\A-[u]_\A\in\Z^d,
\]
where $[v]\in\G$ is the integer part of the vertex $v$ and the vector $[v]_\A\in\Z^d$ is defined by \er{cola}.

For example, for the periodic graph $\cG$ shown in Fig.\ref{ff.0.11}\emph{a} the index of the edge $(v_3+a_2,v_2+a_1)$ is equal to $(1,-1)$, since the integer parts of the vertices $v_3+a_2$ and $v_2+a_1$ are equal to $a_2$ and $a_1$, respectively.

\setlength{\unitlength}{1.0mm}
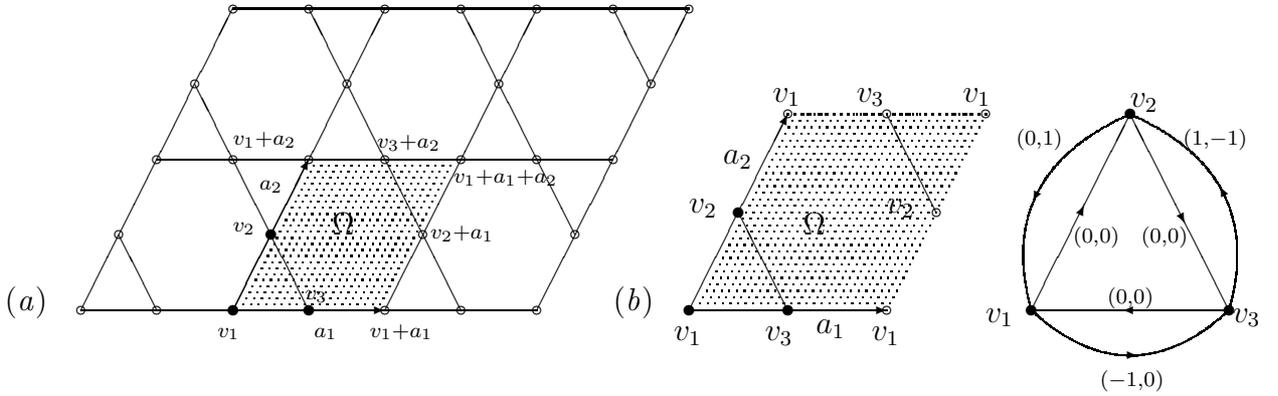
\begin{figure}[h]
\centering
\unitlength 1mm % = 2.845pt
\linethickness{0.4pt}
\ifx\plotpoint\undefined\newsavebox{\plotpoint}\fi % GNUPLOT compatibility

\begin{picture}(160,55)(0,0)
\bezier{30}(81,15)(87.5,28)(94,41)
\bezier{30}(82,15)(88.5,28)(95,41)
\bezier{30}(83,15)(89.5,28)(96,41)
\bezier{30}(84,15)(90.5,28)(97,41)
\bezier{30}(85,15)(91.5,28)(98,41)
\bezier{30}(86,15)(92.5,28)(99,41)
\bezier{30}(87,15)(93.5,28)(100,41)
\bezier{30}(88,15)(94.5,28)(101,41)
\bezier{30}(89,15)(95.5,28)(102,41)
\bezier{30}(90,15)(96.5,28)(103,41)
\bezier{30}(91,15)(97.5,28)(104,41)
\bezier{30}(92,15)(98.5,28)(105,41)
\bezier{30}(93,15)(99.5,28)(106,41)
\bezier{30}(94,15)(100.5,28)(107,41)
\bezier{30}(95,15)(101.5,28)(108,41)
\bezier{30}(96,15)(102.5,28)(109,41)
\bezier{30}(97,15)(103.5,28)(110,41)
\bezier{30}(98,15)(104.5,28)(111,41)
\bezier{30}(99,15)(105.5,28)(112,41)
\bezier{30}(100,15)(106.5,28)(113,41)
\bezier{30}(101,15)(107.5,28)(114,41)
\bezier{30}(102,15)(108.5,28)(115,41)
\bezier{30}(103,15)(109.5,28)(116,41)
\bezier{30}(104,15)(110.5,28)(117,41)
\bezier{30}(105,15)(111.5,28)(118,41)

\bezier{25}(21,15)(26,25)(31,35)
\bezier{25}(22,15)(27,25)(32,35)
\bezier{25}(23,15)(28,25)(33,35)
\bezier{25}(24,15)(29,25)(34,35)
\bezier{25}(25,15)(30,25)(35,35)
\bezier{25}(26,15)(31,25)(36,35)
\bezier{25}(27,15)(32,25)(37,35)
\bezier{25}(28,15)(33,25)(38,35)
\bezier{25}(29,15)(34,25)(39,35)
\bezier{25}(30,15)(35,25)(40,35)
\bezier{25}(31,15)(36,25)(41,35)
\bezier{25}(32,15)(37,25)(42,35)
\bezier{25}(33,15)(38,25)(43,35)
\bezier{25}(34,15)(39,25)(44,35)
\bezier{25}(35,15)(40,25)(45,35)
\bezier{25}(36,15)(41,25)(46,35)
\bezier{25}(37,15)(42,25)(47,35)
\bezier{25}(38,15)(43,25)(48,35)
\bezier{25}(39,15)(44,25)(49,35)

\put(33.0,25){$\Omega$}

\put(95,25){$\Omega$}

\put(0,15){\line(1,0){60.0}}
\put(10,35){\line(1,0){60.0}}
\put(20,55){\line(1,0){60.0}}
\put(0,15){\line(1,2){20.0}}
\put(20,15){\line(1,2){20.0}}
\put(40,15){\line(1,2){20.0}}
\put(60,15){\line(1,2){20.0}}

\put(10,15){\line(-1,2){5.0}}
\put(30,15){\line(-1,2){15.0}}
\put(50,15){\line(-1,2){20.0}}
\put(65,25){\line(-1,2){15.0}}
\put(75,45){\line(-1,2){5.0}}

\put(20,15){\vector(1,0){20.0}}
\put(20,15){\vector(1,2){10.0}}

%\put(20,15.4){\line(1,0){10.0}}
%\put(20,15.2){\line(1,0){10.0}}
%\put(20,14.8){\line(1,0){10.0}}
%\put(20,14.6){\line(1,0){10.0}}

%\put(20,15.6){\line(1,2){5.0}}
%\put(20,15.3){\line(1,2){5.0}}
%\put(20,14.7){\line(1,2){5.0}}
%\put(20,14.4){\line(1,2){5.0}}
%
%\put(30,15.6){\line(-1,2){5.0}}
%\put(30,15.4){\line(-1,2){5.0}}
%\put(30,14.7){\line(-1,2){5.0}}
%\put(30,14.4){\line(-1,2){5.0}}

\put(0,15){\circle{1}}
\put(10,15){\circle{1}}
\put(20,15){\circle*{1.5}}
\put(30,15){\circle*{1.5}}
\put(40,15){\circle{1}}
\put(50,15){\circle{1}}
\put(60,15){\circle{1}}

\put(5,25){\circle{1}}
\put(25,25){\circle*{1.5}}
\put(45,25){\circle{1}}
\put(65,25){\circle{1}}

\put(10,35){\circle{1}}
\put(20,35){\circle{1}}
\put(30,35){\circle{1}}
\put(40,35){\circle{1}}
\put(50,35){\circle{1}}
\put(60,35){\circle{1}}
\put(70,35){\circle{1}}

\put(15,45){\circle{1}}
\put(35,45){\circle{1}}
\put(55,45){\circle{1}}
\put(75,45){\circle{1}}

\put(20,55){\circle{1}}
\put(30,55){\circle{1}}
\put(40,55){\circle{1}}
\put(50,55){\circle{1}}
\put(60,55){\circle{1}}
\put(70,55){\circle{1}}
\put(80,55){\circle{1}}

\put(30.7,11.5){$\scriptstyle a_1$}
\put(23.5,31.0){$\scriptstyle a_2$}
\put(18.0,11.5){$\scriptstyle v_1$}
\put(20.0,37.0){$\scriptstyle v_1+a_2$}
\put(39.0,36.5){$\scriptstyle v_3+a_2$}
\put(38.0,11.5){$\scriptstyle v_1+a_1$}
\put(20.0,25.5){$\scriptstyle v_2$}
\put(46.0,24.5){$\scriptstyle v_2+a_1$}
\put(49.0,32.0){$\scriptstyle v_1+a_1+a_2$}
\put(29.5,16.5){$\scriptstyle v_3$}
\put(-10,15.0){(\emph{a})}

%****************************
%\put(80,15.4){\line(1,0){13.0}}
%\put(80,15.2){\line(1,0){13.0}}
%\put(80,14.8){\line(1,0){13.0}}
%\put(80,14.6){\line(1,0){13.0}}
%
%\put(80,15.3){\line(1,2){6.5}}
%\put(80,15.6){\line(1,2){6.5}}
%\put(80,14.7){\line(1,2){6.5}}
%\put(80,14.4){\line(1,2){6.5}}

\put(80,15){\vector(1,0){26.0}}
\put(80,15){\vector(1,2){13.0}}
\bezier{30}(93,41)(106,41)(119,41)
\bezier{30}(106,15)(112.5,28)(119,41)
\put(80,15){\circle*{1.5}}
\put(93,15){\circle*{1.5}}
\put(106,15){\circle{1}}

\put(93,15){\line(-1,2){6.5}}
%\put(93,15.6){\line(-1,2){6.5}}
%\put(93,15.3){\line(-1,2){6.5}}
%\put(93,14.7){\line(-1,2){6.5}}
%\put(93,14.4){\line(-1,2){6.5}}
\put(112.5,28){\line(-1,2){6.5}}

\put(86.5,28){\circle*{1.5}}
\put(112.5,28){\circle{1}}

\put(93,41){\circle{1}}
\put(106,41){\circle{1}}
\put(119,41){\circle{1}}

\put(96.7,12.0){$a_1$}
\put(84.5,34.5){$a_2$}
\put(78.0,11.0){$v_1$}
\put(91.0,42.7){$v_1$}
\put(102.0,42.7){$v_3$}
\put(104.0,11.0){$v_1$}
\put(80.0,28.0){$v_2$}
\put(106.0,28.0){$v_2$}
\put(116.0,42.7){$v_1$}
\put(90.0,11.0){$v_3$}
\put(70,15.0){(\emph{b})}

%***************************
\put(131.0,27.0){\vector(1,2){1.0}}
\put(144.2,28.5){\vector(1,-2){1.0}}
\put(138,15.0){\vector(-1,0){1.0}}

%\put(125,14.6){\line(1,0){26.0}}
%\put(125,14.8){\line(1,0){26.0}}
\put(125,15){\line(1,0){26.0}}
%\put(125,15.2){\line(1,0){26.0}}
%\put(125,15.4){\line(1,0){26.0}}

%\put(138,14.6){\vector(-1,0){1.0}}
%\put(138,15.4){\vector(-1,0){1.0}}

%\put(125,14.2){\line(1,2){13.0}}
%\put(125,14.6){\line(1,2){13.0}}
\put(125,15){\line(1,2){13.0}}
%\put(125,15.4){\line(1,2){13.0}}
%\put(125,15.8){\line(1,2){13.0}}

%\put(131.5,27.2){\vector(1,2){1.0}}
%\put(131.0,27.5){\vector(1,2){1.0}}
%\put(144.2,28.5){\vector(1,-2){1.0}}

\put(150.8,29.5){\vector(-1,2){1.0}}
\put(126.2,31.5){\vector(-1,-2){1.0}}
\put(138,9.0){\vector(1,0){1.0}}

%\put(151,14.4){\line(-1,2){13.0}}
%\put(151,14.7){\line(-1,2){13.0}}
\put(151,15){\line(-1,2){13.0}}
%\put(151,15.3){\line(-1,2){13.0}}
%\put(151,15.6){\line(-1,2){13.0}}

\put(125,15){\circle*{1.5}}
\put(119.0,13.5){$v_1$}
\put(151,15){\circle*{1.5}}
\put(151.5,13.5){$v_3$}
\put(138,41){\circle*{1.5}}
\put(138.0,42){$v_2$}
\bezier{200}(125,15)(138,3)(151,15)
\bezier{200}(125,15)(120,33)(138,41)
\bezier{200}(151,15)(156,33)(138,41)

\put(135.0,16.0){$\scriptstyle(0,0)$}
\put(130.5,24.0){$\scriptstyle(0,0)$}
\put(139.5,24.0){$\scriptstyle(0,0)$}

\put(134.0,5.0){$\scriptstyle(-1,0)$}
\put(123.0,37.0){$\scriptstyle(0,1)$}
\put(145.0,37.0){$\scriptstyle(1,-1)$}
\end{picture}

\vspace{-0.7cm} \caption{\footnotesize \emph{a}) The Kagome lattice, the fundamental cell $\Omega$ is shaded; \emph{b})~ its fundamental graph with a minimal form $\gm\in\mF(\k)$. The values of $\gm$ are shown near edges.}
\label{ff.0.11}
\end{figure}

On the fundamental graph $G_*$ we introduce the \emph{index form} $\t: \cA_*\ra\Z^d$ by:
\[
\lb{inf}
\t(\bf e_*)=\t(\be) \qq \textrm{ for some $\be\in\cA$ \; such that }  \; \be_*=\gf(\be), \qqq \be_*\in\cA_*,
\]
where $\gf$ is defined by \er{sur}.

The index form $\t$ is uniquely determined by \er{inf}, since
$$
\t(\be+a)=\t(\be),\qqq \forall\, (\be,a)\in\cA \ts\G.
$$

\begin{proposition}
\lb{Ptau}
The index form $\tau:\cA_*\ra\Z^d$ defined by \er{in}, \er{inf} is a 1-form on the fundamental graph $G_*$ and
\[\lb{PaPt}
\Phi_{\tau}=\Phi_{\k},
\]
where $\Phi_{\gb}$ is the flux function of a 1-form $\gb$ defined by \er{mafla}, $\k$ is the coordinate form given by \er{edco}, \er{dco}.
\end{proposition}

This proposition was proved in \cite{KS18} (see Proposition 4.2).

\medskip

\no \textbf{Remark.}  For some periodic graphs index form $\t\in\mF(\k)$ may coincide with a minimal form in $\mF(\k)$.

\subsection{Direct integral decomposition.} Recall that we introduce the Hilbert space $\mH$ by \er{Hisp} and $\lan\cdot\,,\cdot\ran$ denotes the standard inner product in $\R^d$, and $\{a_1,\ldots,a_d\}$ is the basis of the lattice $\G$. We identify the vertices of the fundamental
graph $G_*=(\cV_*,\cE_*)$ with the vertices of the $\G$-periodic graph
$\cG=(\cV,\cE)$ from the fundamental cell $\Omega$. We need Theorem 3.1 from \cite{KS17}.

\begin{theorem}\label{TFD1}
The magnetic Schr\"odinger operator $H_\a=\D_\a+Q$ on $\ell^2(\cV)$ has the following decomposition into a constant fiber direct integral
\[
\lb{raz}
\begin{aligned}
& UH_\a U^{-1}=\int^\oplus_{\T^d}H_{\t,\a}(\vt){d\vt\/(2\pi)^d}
\end{aligned}
\]
for the unitary Gelfand transform $U:\ell^2(\cV)\to\mH$ defined by
\[
\lb{5001}
(Uf)(\vt,v)=\sum\limits_{\mn=(n_1,\ldots,n_d)\in\Z^d}e^{-i\lan \mn,\vt\ran }
f(v+n_1a_1+\ldots+n_da_d), \qqq (\vt,v)\in \T^d\ts\cV_*.
\]
Here the fiber magnetic Schr\"odinger operator $H_{\t,\a}(\vt)$ and the fiber magnetic Laplacian $\D_{\t,\a}(\vt)$ are given by
\[\label{Hvt'}
H_{\t,\a}(\vt)=\D_{\t,\a}(\vt)+Q, \qqq \forall\,\vt\in\T^d,
\]
\[
\label{l2.15''}
 \big(\D_{\t,\a}(\vt)f\big)(v)=\vk_vf(v)-\sum_{\be=(v,\,u)\in\cA_*} e^{i(\a(\be)+\lan\t(\be),\,\vt\ran)}f(u), \qq f\in\ell^2(\cV_*),\qq
 v\in \cV_*,
\]
where $\t$ is the index form defined by \er{in}, \er{inf}.
\end{theorem}

The fiber magnetic Laplacians $\D_{\t,\a}(\cdot)$ given by \er{l2.15''} depend on $\#\cA_*$ indices $\t(\be)$ of the fundamental graph edges $\be\in\cA_*$ and
on $\#\cA_*$ values $\a(\be)$, $\be\in\cA_*$, of the magnetic potentials $\a$. Some of them may be zero, but the number of such zero indices and the number of zero values of the magnetic potential are difficult to control. Moreover, the number of zero indices depends on the choice of the embedding of the periodic graph $\cG$ into the space $\R^d$. Therefore, in Theorem \ref{TDIBB} we construct some gauge transformation and give a representation of the fiber magnetic Laplacians in terms of any 1-forms $(\gb,\ga)\in\mF(\k)\ts\mF(\a)$. In particular, if $\gb$ and $\ga$ are minimal forms, then the fiber operator has the minimal number of coefficients depending on the quasimomentum and the minimal number of coefficients depending on the magnetic form $\ga$. This representation will be used to estimate the Lebesgue measure of the spectrum of the magnetic Schr\"odinger operators on periodic graphs (see Theorem \ref{T1}) and a variation of the spectrum of the Schr\"odinger operators under a perturbation by a magnetic field (see Theorem \ref{Temf}).

\medskip

Let $(\gb,\ga)\in\mF(\k)\ts\mF(\a)$ be two 1-forms on the fundamental graph $G_*$, where $\mF(\mathbf{x})$ is given by \er{vv1f}; $\k$ and $\a$ are the coordinate and magnetic forms defined by \er{edco}, \er{dco} and \er{Mdco}, respectively. We fix a vertex $v_0\in \cV_*$ and introduce the gauge transformation $W_{\gb,\ga}:\mH\to\mH$ by
\[\lb{Uvt}
(W_{\gb,\ga}\,g)(\vt,v)=e^{i(w_{\ga}+\lan w_\gb(v),\vt\ran)}g(\vt,v), \qqq g\in\mH,\qqq  (\vt,v)\in \T^{d}\ts \cV_*,
\]
where the functions $w_\gb:\cV_*\ra\R^d$ and $w_{\ga}:\cV_*\ra\R$ are defined as follows: for any vertex $v\in \cV_*$, we take an oriented path $\textbf{p}=(\be_1,\be_2,\ldots,\be_n)$ on $G_*$ starting at $v_0$ and ending at $v$ and we set
\[\lb{Wvt}
w_\gb(v)=\sum_{s=1}^n\big(\t(\be_s)-\gb(\be_s)\big),\qqq
w_{\ga}(v)=\sum_{s=1}^n\big(\a(\be_s)-\ga(\be_s)\big).
\]

\begin{theorem}
\label{TDIBB}
i) Let $(\gb,\ga)\in\mF(\k)\ts\mF(\a)$, where $\k$ and $\a$ are the coordinate and magnetic forms defined by \er{edco}, \er{dco} and \er{Mdco}, respectively.
Then the magnetic Schr\"odinger operator $H_\a=\D_\a+Q$ on $\ell^2(\cV)$ has the following decomposition into a constant fiber direct integral
\[
\lb{raz1}
\mU_{\gb,\ga} H_\a \mU_{\gb,\ga}^{-1}=\int^\oplus_{\T^d}H_{\gb,\ga}(\vt){d\vt\/(2\pi)^d}\,,
\]
where the unitary operator $\mU_{\gb,\ga}=W_{\gb,\ga}\, U:\ell^2(\cV)\to\mH$, $U$ is the Gelfand transformation given by \er{5001} and $W_{\gb,\ga}$ is the gauge
transformation given by \er{Uvt}, \er{Wvt}.  Here the fiber magnetic Schr\"odinger operator $H_{\gb,\ga}(\vt)$ and the fiber magnetic Laplacian $\D_{\gb,\ga}(\vt)$ are given by
\[
\label{Hvt}
H_{\gb,\ga}(\vt)=\D_{\gb,\ga}(\vt)+Q,\qqq \forall\,\vt\in\T^{d},
\]
\[
\label{l2.13a}
\begin{aligned}
\big(\D_{\gb,\ga}(\vt)f\big)(v)=\vk_vf(v)
-\sum_{\be=(v,u)\in\cA_*}e^{i(\ga(\be)+\lan\gb(\be),\,\vt\ran)}f(u), \qq f\in\ell^2(\cV_*), \qq v\in \cV_*.
\end{aligned}
\]

ii) Let $(\gb',\ga')\in\mF(\k)\ts\mF(\a)$ be two other 1-forms. Then for each $\vt\in\T^{d}$ the operators $\D_{\gb,\ga}(\vt)$ and $\D_{\gb',\ga'}(\vt)$ are unitarily equivalent.

iii) In the identity \er{l2.13a} for the fiber Laplacian $\D_{\gb,\ga}(\cdot)$

$\bu$ the number $\cN_\gb$ of exponents $e^{i\lan\gb(\be),\,\cdot\,\ran}\neq1$, $\be\in\cA_*$, satisfies
\[\lb{dI}
\cN_\gb\geq2\cI, \qqq \textrm{where} \qqq \cI=\textstyle\frac12\,\#\supp \gm=\b(\k),
\]

$\bu$ the number $\cN_\ga$ of exponents $e^{i\ga(\be)}\neq1$, $\be\in\cA_*$, satisfies
\[\lb{dIf}
\cN_\ga\geq2\cI_\a, \qqq \textrm{where} \qqq \cI_\a=\textstyle\frac12\,\#\supp \f=\b(\a),
\]
$(\gm,\f)\in\mF(\k)\ts\mF(\a)$ are two minimal forms, and $\b(\mathbf{x})$ is defined by \er{nbo}.
\end{theorem}

\no{\bf Proof.} \emph{i}) Let $(\gb,\ga)\in\mF(\k)\ts\mF(\a)$. Then, using \er{PaPt}, we have
\[\lb{PbPt}
\Phi_{\gb}=\Phi_{\t},\qqq \Phi_{\ga}=\Phi_{\a}.
\]

Firstly, we show that $w_\gb(v)$ and $w_{\ga}(v)$ defined by \er{Wvt} do not depend on the choice of a path from $v_0$ to $v$. Let $\textbf{p}$ and $\textbf{q}$ be some oriented paths from
$v_0$ to $v$. We consider the cycle $\mathbf{c}=\textbf{p}\ul{\textbf{q}}$, where $\ul{\textbf{q}}$ is the inverse path of $\textbf{q}$. Then we have
\[\lb{com1}
\sum_{\be\in\mathbf{c}}\big(\t(\be)-\gb(\be)\big)=\sum_{\be\in
\textbf{p}}\big(\t(\be)-\gb(\be)\big)-\sum_{\be\in
\textbf{q}}\big(\t(\be)-\gb(\be)\big).
\]
Combining the first identity in \er{PbPt} and \er{com1}, we obtain
$$
\sum_{\be\in \mathbf{p}}\big(\t(\be)-\gb(\be)\big)=\sum_{\be\in
\mathbf{q}}\big(\t(\be)-\gb(\be)\big),
$$
which implies that $w_\gb(v)$ does not depend on the choice of the path from $v_0$ to $v$. Similarly, we can prove that $w_{\ga}(v)$ also does not depend on this choice.

Secondly, we show that the operators $\D_{\t,\a}(\cdot)$ and $\D_{\gb,\ga}(\cdot)$ defined by \er{l2.15''} and \er{l2.13a}, respectively, are unitarily equivalent, by the gauge transformation $W_{\gb,\ga}$ given by \er{Uvt}. Indeed, from \er{Wvt} it follows that
\[\lb{waa}
w_\gb(u)=w_\gb(v)+\t(\be)-\gb(\be), \qq w_{\ga}(u)=w_{\ga}(v)+\a(\be)-\ga(\be), \qq \forall\,\be=(v,u)\in\cA_*.
\]
Using this, \er{l2.15''}, \er{Uvt} and \er{l2.13a}, we have
\[
\begin{aligned}
&\big(\D_{\t,\a}(\vt)f\big)(v)=\vk_vf(v)-\sum_{\be=(v,u)\in\cA_*}
e^{i(\a(\be)+\lan\t(\be),\,\vt\ran)}f(u)\\
&=\vk_vf(v)-\sum_{\be=(v,u)\in\cA_*}e^{-i(w_{\ga}(v)+\lan w_\gb(v),\,\vt\ran)}e^{i(\ga(\be)+\lan\gb(\be),\,\vt\ran)}e^{i(w_{\ga}(u)+\lan w_\gb(u),\,\vt\ran)}f(u)\\
&=e^{-i(w_{\ga}(v)+\lan w_\gb(v),\,\vt\ran)}\Big(\vk_ve^{i(w_{\ga}(v)+\lan w_\gb(v),\,\vt\ran)}f(v)-\hspace{-4mm}\sum_{\be=(v,u)\in\cA_*}
e^{i(\ga(\be)+\lan\gb(\be),\,\vt\ran)}e^{i(w_{\ga}(u)+\lan w_\gb(u),\vt\ran)}f(u)\Big)\\
&=\big(W_{\gb,\ga}^{-1}\D_{\gb,\ga}(\vt)W_{\gb,\ga} f\big)(\vt,v).
\end{aligned}
\]
Thus,
\[\lb{DaDa}
\D_{\t,\a}(\cdot)=W_{\gb,\ga}^{-1}\D_{\gb,\ga}(\cdot)W_{\gb,\ga},
\]
i.e., $\D_{\gb,\ga}(\vt)$ and $\D_{\t,\a}(\vt)$ are unitarily equivalent for each $\vt\in\T^{d}$. Consequently, the fiber Schr\"odinger operator $H_{\t,\a}(\cdot)=\D_{\t,\a}(\cdot)+Q$ is
unitarily equivalent to the operator $H_{\gb,\ga}(\cdot)=\D_{\gb,\ga}(\cdot)+Q$, by the gauge transformation $W_{\gb,\ga}$. This and Theorem \ref{TFD1} give the required statement.

\emph{ii}) This is a direct consequence of the identity \er{DaDa}.

\emph{iii}) Due to the definition of the number $\cN_\gb$, we have $\cN_\gb=\#\supp \gb$. Since $\gm\in\mF(\k)$ is a minimal form,
\[\lb{cNgb}
\cN_\gb=\#\supp \gb\geq \#\supp \gm.
\]
This, \er{dIm} and Theorem \ref{Pphi}.\emph{iii}) give \er{dI}. Similarly, we can prove \er{dIf}. \qq \BBox

\medskip

\no \textbf{Remarks.} 1) The formulas \er{raz1} -- \er{l2.13a} give an infinite number of decompositions into constant fiber direct integrals for the same magnetic Schr\"odinger operator $H_\a$ on the periodic graph $\cG$ (one decomposition for each pair $(\gb,\ga)\in\mF(\k)\ts\mF(\a)$).

2) From Theorem \ref{TDIBB}.\emph{iii}) it follows that the fiber Laplacian $\D_{\gm,\f}(\vt)$ with minimal forms $(\gm,\f)\in\mF(\k)\ts\mF(\a)$ has the minimal number $2\cI$ of coefficients depending on the quasimomentum $\vt$ and the minimal number $2\cI_\a$ of coefficients depending on the 1-form $\f$
among all fiber Laplacians $\D_{\gb,\ga}(\vt)$, $(\gb,\ga)\in\mF(\k)\ts\mF(\a)$.

\medskip

\no{\bf Proof of Theorem \ref{TDImf}.} \emph{i}) The first statement of this item was proved in \cite{KS18} (see Theorem 2.1.\emph{i}). The second one follows from the identities \er{vv1f} and \er{prao}.

\emph{ii}) Since $(\gm,\f)\in\mF(k)\ts\mF(\a)$, the decomposition \er{raz1m} -- \er{l2.13am} is a direct consequence of the decomposition \er{raz1} -- \er{l2.13a} as $\gb=\gm$ and $\ga=\f$.

\emph{iii}) Since $(\gm,\f),(\gm',\f')\in\mF(k)\ts\mF(\a)$, this item follows from Theorem \ref{TDIBB}.\emph{ii}) as $(\gb,\ga)=(\gm,\f)$ and $(\gb',\ga')=(\gm',\f')$.

\emph{iv}) This item is obvious.

\emph{v}) This item was proved in \cite{KS18} (see Theorem 2.1.\emph{v}). \qq \BBox

\begin{corollary}\label{TCo2}
Let $(\gm,\f)\in\mF(\k)\ts\mF(\a)$ be two minimal forms on the fundamental graph $G_*=(\cV_*,\cE_*)$, where $\k$ and $\a$ are the coordinate and magnetic forms defined by \er{edco}, \er{dco} and \er{Mdco}, respectively. Then in the standard orthonormal basis of $\ell^2(\cV_*)=\C^\n$, $\n=\#\cV_*$, the $\n\ts\n$ matrix $\D_{\gm,\f}(\vt)=\{\D_{\gm,\f,uv}(\vt)\}_{u,v\in \cV_*}$ of the fiber Laplacian $\D_{\gm,\f}(\vt)$ defined by \er{l2.13am} is given by
\[\lb{Devt}
\D_{\gm,\f,uv}(\vt)=\left\{
\begin{array}{ll}
\displaystyle\vk_v-\sum\limits_{\be=(v,v)\in\cA_*}\cos(\f(\be)+\lan\gm(\be),\vt\ran), & \textrm{ if } \, u=v\\[20pt]
   \displaystyle-\sum\limits_{\be=(u,v)\in\cA_*}
  e^{-i(\f(\be)+\lan\gm(\be),\vt\ran)},  & \textrm{ if } \, u\sim v,\qq u\neq v\\[20pt]
 \hspace{20mm} 0, \qq & \textrm{ otherwise}\\
\end{array}\right..
\]
Here ${\vk}_v$ is the degree of the vertex $v$.
\end{corollary}

\no {\bf Proof.} Let $(\gh_u)_{u\in\cV_*}$ be the standard orthonormal basis of $\ell^2(\cV_*)$. Substituting the formula \er{l2.13am} into the identity
$$
\D_{\gm,\f,uv}(\vt)=\lan \gh_u,\D_{\gm,\f}(\vt)\gh_v\ran_{\ell^2(\cV_*)}
$$
and using the fact that for each loop $\be=(v,v)\in\cA_*$ there exists a loop $\ul\be=(v,v)\in\cA_*$ and $\gm(\ul\be)=-\gm(\be)$, $\f(\ul\be)=-\f(\be)$,
we obtain \er{Devt}. \qq $\BBox$

\medskip

\no \textbf{Remark.} The number of entries lying on and above the main diagonal of the matrix $\D_{\gm,\f}(\vt)$, defined by \er{Devt}, and depending on $\vt$ can be strictly less than $\cI$, where $\cI=\textstyle\frac12\,\#\supp \gm$, since some edges from $\cE_*\cap\supp \gm$ may be incident to the same pair of vertices. Similarly, the number of entries lying on and above the main diagonal of $\D_{\gm,\f}(\vt)$ and depending on $\f$ can be strictly less than $\cI_\a$, where $\cI_\a=\textstyle\frac12\,\#\supp \f$.

\medskip

\no \textbf{Proof of Proposition \ref{TNNI0}.} \emph{i}) This item was proved in \cite{KS18} (see Proposition 2.2).

\emph{ii})  Let $\mn\in[0,\b]$ be any integer number. We fix some spanning tree $T=(\cV_*,\cE_T)$ of the fundamental graph $G_*=(\cV_*,\cE_*)$. Due to the properties 2) and 3) of spanning trees (see page \pageref{PSTs}), for each $\be\in \cS_T=\cE_*\sm\cE_T$, there exists a unique cycle $\mathbf{c}_\be$ consisting of only $\be$ and edges of $T$ and the set of all these cycles forms a basis of the cycle space $\cC$ of the graph $G_*$. We fix some $\mn$ edges $\be_1,\ldots,\be_\mn\in\cS_T$ and consider a magnetic form $\a:\cA_*\ra(-\pi,\pi]$ on $G_*$ satisfying the following conditions:

$\bu$ $\a(\be_s)$, $s=1,\ldots,\mn$, are rationally independent;

$\bu$ $\a(\be)=0$ for each edge $\be\in\cA_*\sm\{\be_1,\ldots,\be_\mn,\ul\be_1,\ldots,\ul\be_\mn\}$.\\
Due to Theorem \ref{Pphi}.\emph{v}), $\a$ is a minimal form and
$\cI_\a=\frac12\,\#\supp\a=\mn$.

\emph{iii}) Let $M\geq0$, $N>0$ be any integer numbers, and let $G_1$ be a finite connected graph with the Betti number $\b_1=M$. We consider the periodic graph $\cG$ obtained from the $N$-dimensional lattice $\mathbb{L}^N$ by "gluing"\, the graph $G_1=(V,E)$ to each vertex of the lattice (see Fig.\ref{ff.10}). The edge set of the fundamental graph $G_*$ of $\cG$ consists of all edges of the graph $G_1$ and $N$ loops $\be_1,\ldots,\be_N$ of the fundamental graph of $\mathbb{L}^N$. We fix some spanning tree $T=(V,E_T)$ of $G_1$ and consider a magnetic form $\a:\cA_*\ra(-\pi,\pi]$ on $G_*$ satisfying the following conditions:

$\bu$ $\a(\be)=0$ for each edge $\be$ of the tree $T$ and $\a(\be_s)=0$ for each $s=1,\ldots,N$;

$\bu$ $\a(\be)$, $\be\in \cS_T=E\sm E_T$, are rationally independent.

Since $T$ is also a spanning tree of $G_*$, due to Theorem \ref{Pphi}.\emph{v}), $\a$ is a minimal form and
\[\lb{Iamm}
\cI_\a=\textstyle\frac12\,\#\supp\a=\#\cS_T=\b_1=M.
\]
Due to Proposition \ref{TG1}.\emph{i}), the Betti number $\b$ and the invariant $\cI$ of the periodic graph $\cG$ satisfy
$$
\b=M+N, \qqq  \cI=N.
$$
This and \er{Iamm} yield \er{nni00}. 

\emph{iv}) Let $K$ be any nonnegative integer number, and let $G_1$ be a finite connected graph with the Betti number $\b_1=K$. Let $\cG$ be a periodic graph defined in the proof of the previous item. Using Theorem \ref{Pphi}.\emph{i}) -- \emph{ii}) we obtain that the minimal form $\gm\in\mF(\k)$ is unique and $\supp\gm=\{\be_1,\ldots,\be_N,\ul\be_1,\ldots,\ul\be_N\}$, We consider a magnetic form $\a:\cA_*\ra(-\pi,\pi]$ on $G_*$ such that $\supp\a=\supp\gm$. Due to the definition, $\a$ is a minimal form. Since $\supp\a=\supp\gm$, $\cI_{\gm,\f}=N$. Then, using Proposition \ref{TG1}.\emph{i}), we get
$$
\b=\b_1+N=K+\cI_{\gm,\f},
$$
which yields \er{ni00}. \qq \BBox

\medskip

\no \textbf{Proof of Corollary \ref{TCo1}.} Theorem \ref{Pphi}.\emph{i}) gives that the minimal form $\gm$ satisfies $\gm=\wt\gm(\,\cdot\,,T)$ for some $\gm$-minimal spanning tree $T$ of the graph $G$, where $\wt\gm(\,\cdot\,,T):A\ra\R^d$ is defined by \er{cat} as $\mathbf{x}=\gm$. Since $\Phi_{\gm}(\cC)=\Z^d$, we obtain $\gm(A)\ss\Z^d$. We construct a $\Z^d$-periodic graph $\cG$ such that its fundamental graph $G_*$ is isomorphic to $G$ ($G_*\simeq G$), and $\gm$ is the index form on $G_*$. We identify all vertices of the graph $G$ with some distinct points of $[0,1)^d$. Each edge $\be=(u,v)$ of the graph $G$ induces an infinite number of edges
$$
\big\{\big(u+\mn,v+\mn+\gm(\be)\big)\big\}_{\mn\in\Z^d}
$$
of the graph $\cG$ with the edge index $\gm(\be)$. By construction, $\cG$ is a $\Z^d$-periodic graph with the fundamental graph $G_*\simeq G$ and $\gm$ is the index form on $G_*$. Thus, due to Theorem \ref{TFD1}, the operator $\mA(\vt)$, $\vt\in\T^d$, defined by \er{l2.13at} is a fiber operator for the magnetic Laplacian $\D_\f$ with the magnetic potential $\f$ on the graph $\cG$. \qq \BBox

\medskip

\no \textbf{Proof of Corollary \ref{TCo0}.} \emph{i}) Since  $\gm(\be_1),\ldots,\gm(\be_d)\in\Z^d$ are linear independent, there exists
$\vt_0\in\T^d$ satisfying the system of the linear equations
\[\lb{seth}
\f(\be_s)+\lan\gm(\be_s),\vt_0\ran=0,\qqq s=1,\ldots,d.
\]
If we make the change of variables  $\wt\vt=\vt-\vt_0$, then, using \er{seth}, for each $\be\in\cA_*$ we obtain
\begin{multline*}
\f(\be)+\lan\gm(\be),\vt\ran=\f(\be)+\lan\gm(\be),\wt\vt+\vt_0\ran \\
=\left\{
\begin{array}{cl}
0,  & \textrm{ if } \, \be\notin(\supp \gm\cup\supp\phi)\\[6pt]
\lan\gm(\be),\wt\vt\,\ran,  & \textrm{ if } \, \be\in\cB\\[6pt]
\f(\be)+\lan\gm(\be),\vt_0\ran+\lan\gm(\be),\wt\vt\,\ran, \qq & \textrm{ otherwise}\\
\end{array}\right..
\end{multline*}
Thus,
$$
\f(\be)+\lan\gm(\be),\vt\ran=\wt\f(\be)+\lan\gm(\be),\wt\vt\,\ran,\qqq \forall\,\be\in\cA_*,
$$
where $\wt\f$ is defined by \er{wtcat}. This and the decomposition \er{raz1m} -- \er{l2.13am} give \er{raz2m} -- \er{mmfr}.

\emph{ii}) Since $\cB\ss\supp\gm$, we have $\#(\supp \gm\cup\supp\phi)\sm\cB=2(\cI_{\gm,\f}-d)$. But some of $\wt\f(\be)$, $\be\in (\supp \gm\cup\supp\phi)\sm\cB$, may also vanish. This yields
the first statement of the item.

Now let $\cI_{\gm,\f}=d$. Then $\wt\f=0$ and $\D_{\gm,\wt\f}(\vt)=\D_{\gm,0}(\vt)$.
This yields that $H_\a$ is unitarily equivalent to $H_0$. \qq $\BBox$

\section{\lb{Sec4} Proofs of Theorems \ref{T1} and \ref{Temf}}
\setcounter{equation}{0}

In order to prove Theorem \ref{T1} we need the following lemma.

\begin{lemma}\lb{fEst} Let $\vt\in\T^d$, $f\in\ell^2(\cV_*)$,
$\gm\in\mF(\k)$ be a minimal form on the fundamental graph $G_*=(\cV_*,\cE_*)$, where $\mF(\k)$ is given by \er{vv1f} as $\mathbf{x}=\k$; $\k$ is the coordinate form defined by \er{edco}, \er{dco}. Let $\a$ be the magnetic form on $G_*$ defined by \er{Mdco}. Then

i) The fiber magnetic Laplacian $\D_{\gm,\a}(\vt)$ given by \er{l2.13a} as $(\gb,\ga)=(\gm,\a)$, has the following representation:
\[\lb{fDmDb}
\D_{\gm,\a}(\vt)=\D_\a^0+\wt\D_{\gm,\a}(\vt),
\]
where $\D_\a^0$ is the magnetic Laplacian with the magnetic vector potential $\a$ on the finite graph $G_\gm^0=(\cV_*,\cE_*\sm\cE_\gm)$, $\cE_\gm=\cE_*\cap\supp \gm$;
$\wt\D_{\gm,\a}(\vt)$ is the magnetic Laplacian with the magnetic vector potential $\a(\be)+\lan\gm(\be),\vt\ran$, $\be\in\cA_\gm$, on the finite graph
$G_\gm=(\cV_*,\cE_\gm)$:
\[\label{fDbt}
\big(\wt\D_{\gm,\a}(\vt)f\big)(v)=\sum_{\be=(v,\,u)\in\cA_\gm}
\big(f(v)-e^{i(\a(\be)+\lan\gm(\be),\,\vt\ran)}f(u)\big), \qq
v\in \cV_*,
\]
$\cA_\gm=\{\be \mid \be\in\cE_\gm \; \textrm{or} \;\; \ul\be\in\cE_\gm\}$.

ii) The quadratic form of the magnetic Laplacian
$\wt\D_{\gm,\a}(\vt)$ is given by
\[\lb{qflom}
\lan \wt\D_{\gm,\a}(\vt)
f,f\ran_{\ell^2(\cV_*)}={1\/2}\sum_{\be=(v,u)\in\cA_\gm}
\big|f(v)-e^{i(\a(\be)+\lan\gm(\be),\,\vt\ran)}f(u)\big|^2.
\]

iii) The magnetic Laplacian $\wt\D_{\gm,\a}(\vt)$ given by \er{fDbt} satisfies
\[\lb{0DaB}
0\le\wt\D_{\gm,\a}(\vt)\le2B_\gm,
\]
where
\[\lb{Bfv}
(B_\gm f)(v)=\vk_v^{\gm}f(v), \qqq v\in \cV_*,
\]
$\vk_v^{\gm}$ is the degree of the vertex $v\in \cV_*$ on the graph $G_\gm=(\cV_*,\cE_\gm)$, and
\[\lb{fsdbe}
\s\big(\wt\D_{\gm,\a}(\vt)\big)\ss[0,2\vk_+^{\gm}], \qqq \vk^{\gm}_+=\max_{v\in \cV_*}\vk_v^{\gm}.
\]
\end{lemma}

\no \textbf{Proof.} \emph{i}) Let $f\in\ell^2(\cV_*)$, $v\in \cV_*$. Then using \er{l2.13a} for each $\vt\in\T^d$ we obtain
$$
\begin{aligned}
&\big(\D_{\gm,\a}(\vt)f\big)(v)=\sum_{\be=(v,\,u)\in\cA_*}
\big(f(v)-e^{i(\a(\be)+\lan\gm(\be),\,\vt\ran)}f(u)\big)\\&=
\sum_{\be=(v,\,u)\in\cA_*\sm\cA_\gm}\big(f(v)-e^{i\a(\be)}f(u)\big)+
\sum_{\be=(v,\,u)\in\cA_\gm}\big(f(v)-e^{i(\a(\be)+\lan\gm(\be),\,\vt\ran)}f(u)\big)=
\D_\a^0+\wt\D_{\gm,\a}(\vt).
\end{aligned}
$$
The operator $\D_\a^0$ is the magnetic Laplacian with the magnetic potential $\a(\be)$, $\be\in\cA_*\sm\cA_\gm$, on the finite graph $G_\gm^0=(\cV_*,\cE_*\sm\cE_\gm)$, and for each $\vt\in\T^d$ the operator $\wt\D_{\gm,\a}(\vt)$ is the magnetic Laplacian with the magnetic potential $\a(\be)+\lan\gm(\be),\vt\ran$, $\be\in\cA_\gm$, on the finite graph $G_\gm=(\cV_*,\cE_\gm)$.

\emph{ii}) This item was proved in \cite{KS17} (see Theorem 6.4.iv).

\emph{iii}) We show that $2B_\gm-\wt\D_{\gm,\a}(\vt)\ge0$ for each $\vt\in\T^d$. Let $f\in\ell^2(\cV_*)$. Then, using the identities \er{qflom}, \er{Bfv}, we obtain
\[\lb{2BD0}
\begin{aligned}
&\textstyle\lan\big(2B_\gm-\wt\D_{\gm,\a}(\vt)\big)f,f\ran_{\ell^2(\cV_*)}=
2\sum\limits_{v\in \cV_*}\vk_v^\gm|f(v)|^2-\lan\wt\D_{\gm,\a}(\vt)f,f\ran_{\ell^2(\cV_*)}\\
&\textstyle=\sum\limits_{(v,\,u)\in\cA_\gm}\big(|f(v)|^2+
|f(u)|^2\big)-{1\/2}\sum\limits_{\be=(v,\,u)\in\cA_\gm}
\big|f(v)-e^{i(\a(\be)+\lan\gm(\be),\,\vt\ran)}f(u)\big|^2\\
&\textstyle={1\/2}\sum\limits_{\be=(v,\,u)\in\cA_\gm}
\big|f(v)+e^{i(\a(\be)+\lan\gm(\be),\,\vt\ran)}f(u)\big|^2\geq0.
\end{aligned}
\]
For each\, $\vt\in\T^d$ the spectrum of the magnetic operator $\wt\D_{\gm,\a}(\vt)$  satisfies the condition \er{fsdbe} (see, e.g., \cite{HS99a}), which, in particularly, yields that $\wt\D_{\gm,\a}(\vt)\geq0$. \qq \BBox

\medskip

\no {\bf Proof of Theorem \ref{T1}.} \emph{i}) Using \er{fDmDb}, we rewrite the fiber magnetic Schr\"odinger operator $H_{\gm,\a}(\vt)$,  $\vt\in\T^d$, defined by \er{Hvt}, \er{l2.13a} as $(\gb,\ga)=(\gm,\a)$, in the form:
\[
\label{eq.1}
H_{\gm,\a}(\vt)=H_\a^0+\wt\D_{\gm,\a}(\vt),
\]
where $H_\a^0=\D_\a^0+Q$ is the magnetic Schr\"odinger operator on the finite graph $G_\gm^0=(\cV_*,\cE_*\sm\cE_\gm)$, $\wt\D_{\gm,\a}(\vt)$ is the magnetic Laplacian on the finite graph $G_\gm=(\cV_*,\cE_\gm)$ defined by \er{fDbt}. Lemma \ref{fEst}.\emph{iii} gives that the spectrum $\s\big(\wt\D_{\gm,\a}(\vt)\big)\ss[0,2\vk^\gm_+]$, where $\vk^\gm_+$ is defined in \er{fesbp1}. Then each eigenvalue $\l_{\a,n}(\vt)$, $n\in\N_\n$, of $H_{\gm,\a}(\vt)$ satisfies $\m_{\a,n}\le\l_{\a,n}(\vt)\le\m_{\a,n}+2\vk^\gm_+$, which yields
\er{fesbp1}.

\emph{ii}) From \er{eq.1} and \er{0DaB} we obtain
$$
H_\a^0\le H_{\gm,\a}(\vt)\le H_\a^0+2B_\gm,
$$
where $B_\gm$ is defined in \er{Bfv}. This yields
$$
\l_{n}(H_\a^0)\leq\l_{\a,n}^-\le \l_{\a,n}(\vt)\leq\l_{\a,n}^+\le \l_{n}(H_\a^0+2B_\gm),
\qqq \forall\,(n,\vt)\in\N_\n\ts\T^d.
$$
Then, using the first identity in \er{dIm},
\begin{multline*}
\big|\s(H_\a)\big|\le\sum_{n=1}^{\nu}(\l_{\a,n}^+-\l_{\a,n}^-)\leq
\sum_{n=1}^\nu\big(\l_{n}(H_\a^0+2B_\gm)-\l_{n}(H_\a^0)\big)\\=2\Tr B_\gm=2\sum_{v\in\cV_*}\vk^\gm_v=
4\#\cE_\gm=2\#\supp \gm=4\cI.
\end{multline*}

\emph{iii}) This item follows from Proposition \ref{TG1}.\emph{ii}).
\qq $\BBox$

%*************************************************************
\setlength{\unitlength}{1.0mm}
\begin{figure}[h]
\centering
\unitlength 1mm % = 2.845pt
\linethickness{0.4pt}
\ifx\plotpoint\undefined\newsavebox{\plotpoint}\fi % GNUPLOT compatibility
\begin{picture}(120,50)(0,0)
\put(0,40){$\cG$}
\put(0,10){\emph{a})}
\put(20,20){\vector(1,0){20.00}}
\put(20,20){\vector(0,1){20.00}}
\multiput(21,40)(4,0){5}{\line(1,0){2}}
\multiput(40,21)(0,4){5}{\line(0,1){2}}
\put(38.5,18.0){$\scriptstyle a_1$}
\put(16.5,39.0){$\scriptstyle a_2$}
%\put(18,18){$\scriptstyle O$}
\put(35,25){$\Omega$}

\bezier{40}(20.5,20)(20.5,30)(20.5,40)
\bezier{40}(21.0,20)(21.0,30)(21.0,40)
\bezier{40}(21.5,20)(21.5,30)(21.5,40)
\bezier{40}(22.0,20)(22.0,30)(22.0,40)
\bezier{40}(22.5,20)(22.5,30)(22.5,40)
\bezier{40}(23.0,20)(23.0,30)(23.0,40)
\bezier{40}(23.5,20)(23.5,30)(23.5,40)
\bezier{40}(24.0,20)(24.0,30)(24.0,40)
\bezier{40}(24.5,20)(24.5,30)(24.5,40)
\bezier{40}(25.0,20)(25.0,30)(25.0,40)
\bezier{40}(25.5,20)(25.5,30)(25.5,40)
\bezier{40}(26.0,20)(26.0,30)(26.0,40)
\bezier{40}(26.5,20)(26.5,30)(26.5,40)
\bezier{40}(27.0,20)(27.0,30)(27.0,40)
\bezier{40}(27.5,20)(27.5,30)(27.5,40)
\bezier{40}(28.0,20)(28.0,30)(28.0,40)
\bezier{40}(28.5,20)(28.5,30)(28.5,40)
\bezier{40}(29.0,20)(29.0,30)(29.0,40)
\bezier{40}(29.5,20)(29.5,30)(29.5,40)
\bezier{40}(30.0,20)(30.0,30)(30.0,40)

\bezier{40}(30.5,20)(30.5,30)(30.5,40)
\bezier{40}(31.0,20)(31.0,30)(31.0,40)
\bezier{40}(31.5,20)(31.5,30)(31.5,40)
\bezier{40}(32.0,20)(32.0,30)(32.0,40)
\bezier{40}(32.5,20)(32.5,30)(32.5,40)
\bezier{40}(33.0,20)(33.0,30)(33.0,40)
\bezier{40}(33.5,20)(33.5,30)(33.5,40)
\bezier{40}(34.0,20)(34.0,30)(34.0,40)
\bezier{40}(34.5,20)(34.5,30)(34.5,40)
\bezier{40}(35.0,20)(35.0,30)(35.0,40)
\bezier{40}(35.5,20)(35.5,30)(35.5,40)
\bezier{40}(36.0,20)(36.0,30)(36.0,40)
\bezier{40}(36.5,20)(36.5,30)(36.5,40)
\bezier{40}(37.0,20)(37.0,30)(37.0,40)
\bezier{40}(37.5,20)(37.5,30)(37.5,40)
\bezier{40}(38.0,20)(38.0,30)(38.0,40)
\bezier{40}(38.5,20)(38.5,30)(38.5,40)
\bezier{40}(39.0,20)(39.0,30)(39.0,40)
\bezier{40}(39.5,20)(39.5,30)(39.5,40)
\bezier{40}(40.0,20)(40.0,30)(40.0,40)

\put(10,10){\line(1,0){40.00}}
\put(10,30){\line(1,0){40.00}}
\put(30,10){\line(0,1){40.00}}
\put(10,50){\line(1,0){40.00}}
\put(10,10){\line(0,1){40.00}}
\put(50,10){\line(0,1){40.00}}

\put(10,10){\circle{1}}
\put(30,10){\circle{1}}
\put(50,10){\circle{1}}

\put(10,30){\circle{1}}
\put(30,30){\circle*{1}}
\put(50,30){\circle{1}}

\put(10,50){\circle{1}}
\put(30,50){\circle{1}}
\put(50,50){\circle{1}}
\put(30,30){\line(1,1){4.00}}
\put(34,34){\line(1,0){4.00}}
\put(34,34){\line(0,1){4.00}}
\put(34,38){\line(1,-1){4.00}}
\put(34,34){\circle*{1}}
\put(34,38){\circle*{1}}
\put(38,34){\circle*{1}}

\put(30,30){\line(-1,1){4.00}}
\put(26,34){\line(-1,0){4.00}}
\put(26,34){\line(0,1){4.00}}
\put(26,38){\line(-1,-1){4.00}}
\put(26,34){\circle*{1}}
\put(26,38){\circle*{1}}
\put(22,34){\circle*{1}}

\put(30,30){\line(-1,-2){3.00}}
\put(30,30){\line(1,-2){3.00}}
\put(27,24){\line(1,0){6.00}}
\put(27,24){\circle*{1}}
\put(33,24){\circle*{1}}
%************************

\put(50,30){\line(-1,1){4.00}}
\put(46,34){\line(-1,0){4.00}}
\put(46,34){\line(0,1){4.00}}
\put(46,38){\line(-1,-1){4.00}}
\put(46,34){\circle{1}}
\put(46,38){\circle{1}}
\put(42,34){\circle{1}}

\put(50,30){\line(-1,-2){3.00}}
\put(47,24){\line(1,0){3.00}}
\put(47,24){\circle{1}}
%************************
\put(10,30){\line(1,1){4.00}}
\put(14,34){\line(1,0){4.00}}
\put(14,34){\line(0,1){4.00}}
\put(14,38){\line(1,-1){4.00}}
\put(14,34){\circle{1}}
\put(14,38){\circle{1}}
\put(18,34){\circle{1}}

\put(10,30){\line(1,-2){3.00}}
\put(10,24){\line(1,0){3.00}}
\put(13,24){\circle{1}}
%***************************
\put(30,10){\line(1,1){4.00}}
\put(34,14){\line(1,0){4.00}}
\put(34,14){\line(0,1){4.00}}
\put(34,18){\line(1,-1){4.00}}
\put(34,14){\circle{1}}
\put(34,18){\circle{1}}
\put(38,14){\circle{1}}

\put(30,10){\line(-1,1){4.00}}
\put(26,14){\line(-1,0){4.00}}
\put(26,14){\line(0,1){4.00}}
\put(26,18){\line(-1,-1){4.00}}
\put(26,14){\circle{1}}
\put(26,18){\circle{1}}
\put(22,14){\circle{1}}

%************************
\put(50,10){\line(-1,1){4.00}}
\put(46,14){\line(-1,0){4.00}}
\put(46,14){\line(0,1){4.00}}
\put(46,18){\line(-1,-1){4.00}}
\put(46,14){\circle{1}}
\put(46,18){\circle{1}}
\put(42,14){\circle{1}}

%************************
\put(10,10){\line(1,1){4.00}}
\put(14,14){\line(1,0){4.00}}
\put(14,14){\line(0,1){4.00}}
\put(14,18){\line(1,-1){4.00}}
\put(14,14){\circle{1}}
\put(14,18){\circle{1}}
\put(18,14){\circle{1}}
%***********************
\put(30,50){\line(-1,-2){3.00}}
\put(30,50){\line(1,-2){3.00}}
\put(27,44){\line(1,0){6.00}}
\put(27,44){\circle{1}}
\put(33,44){\circle{1}}
%************************
\put(50,50){\line(-1,-2){3.00}}
\put(47,44){\line(1,0){3.00}}
\put(47,44){\circle{1}}
%************************
\put(10,50){\line(1,-2){3.00}}
\put(10,44){\line(1,0){3.00}}
\put(13,44){\circle{1}}

%**********************************
\put(58,10){\emph{b})}
%*************************************
\put(100,30){\circle*{1}}

\put(100,30){\line(1,1){4.00}}
\put(104,34){\line(1,0){4.00}}
\put(104,34){\line(0,1){4.00}}

\put(104,38){\line(1,-1){4.00}}
\put(104,34){\circle*{1}}
\put(104,38){\circle*{1}}
\put(108,34){\circle*{1}}

\put(100,30){\line(-1,1){4.00}}
\put(96,34){\line(-1,0){4.00}}
\put(96,34){\line(0,1){4.00}}

\put(96,38){\line(-1,-1){4.00}}
\put(96,34){\circle*{1}}
\put(96,38){\circle*{1}}
\put(92,34){\circle*{1}}

\put(100,30){\line(-1,-2){3.00}}
\put(100,30){\line(1,-2){3.00}}

\put(97,24){\circle*{1}}
\put(103,24){\circle*{1}}
\put(97,24){\line(1,0){6.00}}

\bezier{300}(100,30)(105,43)(100,44)
\bezier{300}(100,30)(95,43)(100,44)

\bezier{300}(100,30)(113,35)(114,30)
\bezier{300}(100,30)(113,25)(114,30)

\put(87,40){$G_*$}
\put(97,45){$\scriptstyle(0,1)$}
\put(106,24.5){$\scriptstyle(1,0)$}
\put(85,10){\emph{c}) }
%********************************

\put(57,40){$G_1$}
\put(70,30){\line(1,1){4.00}}
\put(74,34){\line(1,0){4.00}}
\put(74,34){\line(0,1){4.00}}
\put(74,38){\line(1,-1){4.00}}
\put(74,34){\circle*{1}}
\put(74,38){\circle*{1}}
\put(78,34){\circle*{1}}

\put(70,30){\line(-1,1){4.00}}
\put(66,34){\line(-1,0){4.00}}
\put(66,34){\line(0,1){4.00}}
\put(66,38){\line(-1,-1){4.00}}
\put(66,34){\circle*{1}}
\put(66,38){\circle*{1}}
\put(62,34){\circle*{1}}

\put(70,30){\line(-1,-2){3.00}}
\put(70,30){\line(1,-2){3.00}}
\put(67,24){\line(1,0){6.00}}
\put(67,24){\circle*{1}}
\put(73,24){\circle*{1}}

%*********************************
%\put(64,10){\emph{d})}
%
%\put(80,10){\line(1,0){50.00}}
%\put(80,9){\line(0,1){2.00}}
%\put(95,9){\line(0,1){2.00}}
%\put(84,9){\line(0,1){2.00}}
%\put(130,9){\line(0,1){2.00}}
%\put(85,10){\circle*{1}}
%
%\put(80,9.8){\line(1,0){4.00}}
%\put(80,10.2){\line(1,0){4.00}}
%
%\put(95,9.8){\line(1,0){35.00}}
%\put(95,10.2){\line(1,0){35.00}}
%
%%\put(110,12){$\s_3$}
%%\put(80,11.5){$\scriptstyle\s_1$}
%\put(84.3,11.5){$\scriptstyle1$}
%\put(79,6){$\scriptstyle0$}
%\put(95,6){$\scriptstyle3$}
%\put(123.0,5.5){$\scriptstyle\frac{11+\sqrt{89}}2$}
%\put(82.0,5.5){$\scriptstyle\frac{11-\sqrt{89}}2$}
\end{picture}
\vspace{-0.5cm}
\caption{\footnotesize  \emph{a}) $\Z^2$-periodic graph $\cG$, $a_1,a_2$ are the basis of the lattice $\Z^2$; \quad \emph{b}) a finite decoration $G_1$;\quad \emph{c}) the fundamental graph $G_*$, the values $\gm(\be)$, $\be\in\supp\gm$, of the minimal form $\gm\in\mF(\k)$ are shown near the edges.}
\label{ff.10}
\end{figure}
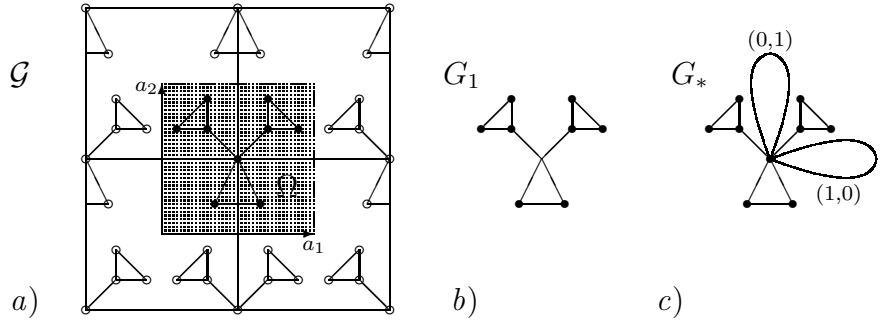

\begin{proposition}\lb{TG1}
Let $\cG$ be a periodic graph obtained from the
$d$-dimensional lattice by "gluing"\, the same finite connected graph $G_1$ with the Betti number $\b_1$ to each vertex of the lattice (for $d=2$ see Fig.\ref{ff.10}a). Then the following statements hold true.

i) The Betti number $\b$ and the invariant $\cI$ given in \er{dIm} for the periodic graph $\cG$ satisfy
\[\lb{deco}
\b=\b_1+d, \qqq \cI=d.
\]

ii) Let, in addition, $G_1$ be a tree. Then the Lebesgue measure of the spectrum of the magnetic Schr\"odinger operator $H_\a$ on $\cG$ satisfies
\[
\lb{sp2}
\textstyle|\s(H_\a)|=4d=4\cI
\]
and the estimate \er{eq.7'} becomes an identity.
\end{proposition}

\no {\bf Proof.} \emph{i}) This item was proved in \cite{KS18} (see Proposition 5.2.i).

\emph{ii}) Since $G_1$ is a tree, we have $\b_1=0$ and, using \er{deco}, $\b=d$.
The last identity gives that the magnetic Schr\"odinger operator $H_{\a}$ is unitarily equivalent to the Schr\"odinger operator $H_0$ without a magnetic field (see, e.g., Corollary 2.2 in \cite{KS17}). The identity $\textstyle|\s(H_0)|=4d$ were proved in \cite{KS14} (see Proposition 7.2). Since $d=\cI$, for the considered graph $\cG$ the estimate \er{eq.7'} becomes an identity. \qq \BBox

\medskip

\no {\bf Proof of Theorem \ref{Temf}.} We rewrite the fiber magnetic Laplacian $\D_{\gm,\wt\f}(\vt)$ defined by \er{mmfr} in the form
\[
\label{l2.13am2}
\D_{\gm,\wt\f}(\vt)=\D_{\gm,0}(\vt)+X_{\gm,\wt\f}(\vt),
\]
where the operator $X_{\gm,\wt\f}(\vt)$, $\vt\in\T^d$, is given by
\[\lb{Xmf}
\big(X_{\gm,\wt\f}(\vt)f\big)(v)=\sum_{\be=(v,u)\in\supp\wt\f}e^{i\lan\gm(\be),\vt\ran}
\big(1-e^{i\wt\phi(\be)}\big)f(u), \qq f\in\ell^2(\cV_*), \qq v\in \cV_*.
\]
Let $\l_{\wt\f,1}(\vt)\le\l_{\wt\f,2}(\vt)\le\ldots\le \l_{\wt\f,\n }(\vt)$ be the eigenvalues of $H_{\gm,\wt\f}(\vt)=\D_{\gm,\wt\f}(\vt)+Q$. Due to \er{l2.13am2},
we have
\[
H_{\gm,\wt\f}(\vt)=H_{\gm,0}(\vt)+X_{\gm,\wt\f}(\vt).
\]
This yields
\[
\lb{es1} \l_{0,n}(\vt)+\L_1\leq\l_{\wt\f,n}(\vt)\leq\l_{0,n}(\vt)+\L_\n, \qqq n\in\N_\n,
\]
where $\L_1$, $\L_\n$ are defined by \er{emf3'}. From \er{es1} and the identities
$$
\l_{\a,n}^+=\max_{\vt\in\T^d}\l_{\wt\f,n}(\vt),\qqq \l_{\a,n}^-=\min_{\vt\in\T^d}\l_{\wt\f,n}(\vt)
$$
we deduce that
$\L_1\le \l_{\a,n}^\pm-\l_{0,n}^{\pm}\le\L_\n$ and
\[
\lb{ees3} |\s_{\a,n}|=\l_{\a,n}^+-\l_{\a,n}^-\le
(\l_{0,n}^{+}-\l_{0,n}^{-})
+(\L_\n-\L_1)=|\s_{0,n}|+(\L_\n-\L_1).
\]
Similar arguments give
\[
\lb{ees4}
|\s_{\a,n}|=\l_{\a,n}^+-\l_{\a,n}^-\geq
(\l_{0,n}^{+}-\l_{0,n}^{-})
-(\L_\n-\L_1)=|\s_{0,n}|-(\L_\n-\L_1).
\]
Combining \er{ees3} and \er{ees4}, we obtain \er{emf2}.

We estimate $\L_1$ and $\L_\n$. The standard estimate and the formula \er{Xmf} yield
$$
\max\{|\L_1|,|\L_\n|\}\le \max_{\vt\in\T^d}\|X_{\gm,\wt\f}(\vt)\|,
$$
$$
\begin{aligned}
&\textstyle\|X_{\gm,\wt\f}(\vt)\|\le
\max\limits_{v\in\cV_*} \sum\limits_{u\in\cV_*}
\Big|\sum\limits_{\be=(v,u)\in\supp\wt\f}e^{i\lan\gm(\be),\vt\ran}
\big(1-e^{i\wt\phi(\be)}\big)\Big|\\
&\textstyle\le\max\limits_{v\in\cV_*} \sum\limits_{u\in\cV_*}\sum\limits_{\be=(v,u)\in\supp\wt\f}
|1-e^{i\wt\phi(\be)}|=\max\limits_{v\in\cV_*}\sum\limits_{\be=(v,u)\in\supp\wt\f}
2\,\big|\sin\frac{\wt\phi(\be)}2\big|.
\end{aligned}
$$
Then we deduce that
\[
\lb{es3.1}
\textstyle\max\limits\{|\L_1|,|\L_\n|\}\le\max\limits_{\vt\in \T^d}\|X_{\gm,\wt\f}(\vt)\|\le
\max\limits_{v\in\cV_*}\sum\limits_{\be=(v,u)\in\supp\wt\f}
2\,\big|\sin\frac{\wt\phi(\be)}2\big|=C_{\wt\f}\,,
\]
where $C_{\wt\f}$ is defined in \er{es3}. The second identity in \er{es30} is a simple consequence of the first one.
$\BBox$

\medskip

\footnotesize
\textbf{Acknowledgments. \lb{Sec8}}  Our study was
supported by the RSF grant  No. 18-11-00032.

%*****************************************************

\end{document}